%
% Projectivity of Hopf algebras over subalgebras with semilocal central localizations
%	S. Skryabin
%
%	Plain TeX + AMS fonts ( amssym.def )
%
%	PARAMETERS
%
\hsize=5in
\baselineskip=12pt
\vsize=20.5cm
\parindent=.5cm
\predisplaypenalty=0
\hfuzz=2pt
\frenchspacing
%
%	Additional Fonts
%
\input amssym.def
\def\titlefonts{\baselineskip=1.44\baselineskip
	\font\titlef=cmbx12
	\titlef
	}
\font\ninerm=cmr9
\font\ninebf=cmbx9
\font\ninei=cmmi9
\skewchar\ninei='177
\font\ninesy=cmsy9
\skewchar\ninesy='60
\font\nineit=cmti9
\def\reffonts{\baselineskip=0.9\baselineskip
	\textfont0=\ninerm
	\def\rm{\fam0\ninerm}%
	\textfont1=\ninei
	\textfont2=\ninesy
	\textfont\bffam=\ninebf
	\def\bf{\fam\bffam\ninebf}%
	\def\it{\nineit}%
	}
%
%	Formatting
%
\def\frontmatter{\vbox{}\vskip1cm\bgroup
	\leftskip=0pt plus1fil\rightskip=0pt plus1fil
	\parindent=0pt
	\parfillskip=0pt
	\pretolerance=10000
	}
\def\endfrontmatter{\egroup\bigskip}
\def\title#1{{\titlefonts#1\par}}
\def\author#1{\bigskip#1\par}
\def\address#1{\bigskip{\reffonts\it#1}}
\def\email#1{\bigskip{\reffonts{\it E-mail: }\rm#1}}
\def\thanks#1{\footnote{}{\reffonts\rm\noindent#1\hfil}}
\def\section#1\par{\ifdim\lastskip<\bigskipamount\removelastskip\fi
	\penalty-250\bigskip
	\vbox{\leftskip=0pt plus1fil\rightskip=0pt plus1fil
	\parindent=0pt
	\parfillskip=0pt
	\pretolerance=10000{\bf#1}}\nobreak\medskip
	}
\def\proclaim#1. {\medbreak\bgroup{\noindent\bf#1.}\ \it}
\def\endproclaim{\egroup
	\ifdim\lastskip<\medskipamount\removelastskip\medskip\fi}
\newdimen\itemsize
\def\setitemsize#1 {{\setbox0\hbox{#1\ }
	\global\itemsize=\wd0}}
\def\item#1 #2\par{\ifdim\lastskip<\smallskipamount\removelastskip\smallskip\fi
	{\leftskip=\itemsize
	\noindent\hskip-\leftskip
	\hbox to\leftskip{\hfil\rm#1\ }#2\par}\smallskip}
\def\Proof#1. {\ifdim\lastskip<\medskipamount\removelastskip\medskip\fi
	{\noindent\it Proof\if\space#1\space\else\ \fi#1.}\ }
\def\endproof{\quad\hbox{$\square$}\medskip}
\def\Remark. {\ifdim\lastskip<\medskipamount\removelastskip\medskip\fi
	{\noindent\bf Remark.}\quad}
\def\endremark{\medskip}
%
%		Citations
%
\newcount\citation
\newtoks\citetoks
\def\citedef#1\endcitedef{\citetoks={#1\endcitedef}}
\def\endcitedef#1\endcitedef{}
\def\citenum#1{\citation=0\def\curcite{#1}%
	\expandafter\checkendcite\the\citetoks}
\def\checkendcite#1{\ifx\endcitedef#1?\else
	\expandafter\lookcite\expandafter#1\fi}
\def\lookcite#1 {\advance\citation by1\def\auxcite{#1}%
	\ifx\auxcite\curcite\the\citation\expandafter\endcitedef\else
	\expandafter\checkendcite\fi}
\def\cite#1{\makecite#1,\cite}
\def\makecite#1,#2{[\citenum{#1}\ifx\cite#2]\else\expandafter\clearcite\expandafter#2\fi}
\def\clearcite#1,\cite{, #1]}
\citedef
Az51
Bass
Bj73
Bou
Br
Camps93
Ch90
Cor88
Cur53
Con94
Doi83
Et04
Hoff
Kh
Lam
Ma91
Ma92
Ma94
Mo
Mo95
Mo99
Ni89
Pro67
Rad77
Rob81
Scha00
Schn93
Sk07
Sk-Oy07
Sw
Tak71
Tak79
Wu03
\endcitedef
%
%		Reference section
%
\def\references{\section References\par
	\bgroup
	\parindent=0pt
	\reffonts
	\rm
	\frenchspacing
	\setbox0\hbox{99. }\leftskip=\wd0
	}
\def\endreferences{\egroup}
\newtoks\nextauth
\newif\iffirstauth
\def\checkendauth#1{\ifx\endauth#1%
		\iffirstauth\the\nextauth
		\else{} and \the\nextauth\fi,
	\else\iffirstauth\the\nextauth\firstauthfalse
		\else, \the\nextauth\fi
		\expandafter\auth\expandafter#1\fi
	}
\def\auth#1,#2;{\nextauth={#1 #2}\checkendauth}
\newif\ifinbook
\newif\ifbookref
\def\nextref#1 {\par\hskip-\leftskip
	\hbox to\leftskip{\hfil\citenum{#1}.\ }%
	\initnextref}
\def\initnextref{\bookreffalse\inbookfalse\firstauthtrue\ignorespaces}
\def\paper#1{{\it#1},}
\def\In#1{\inbooktrue in ``#1",}
\def\book#1{\bookreftrue{\it#1},}
\def\journal#1{#1\ifinbook,\fi}
\def\bookseries#1{#1,}
\def\Vol#1{\ifbookref Vol. #1,\else\ifinbook Vol. #1,\else{\bf#1}\fi\fi}

\def\publisher#1{#1,}
\def\Year#1{\ifbookref #1.\else\ifinbook #1,\else(#1)\fi\fi}
\def\Pages#1{\makepages#1.}
\long\def\makepages#1-#2.#3{\ifinbook pp. \fi#1--#2\if\par#3.\fi#3}

%
%	Math Stuff
%
\newsymbol\square 1003
\let\ot\otimes
\let\sbs\subset
\let\sps\supset
\newsymbol\rightsquigarrow 1320
\newsymbol\smallsetminus 2272
\let\setm\smallsetminus
\newsymbol\varnothing 203F
\def\ev{{\rm ev}}
\def\Hom{\mathop{\rm Hom}\nolimits}
\def\id{{\rm id}}
\def\inf{\mathop{\rm inf}\nolimits}
\def\Jac{\mathop{\rm Jac}\nolimits}
\def\Ker{\mathop{\rm Ker}}
\def\limdir{\mathop{\vtop{\offinterlineskip\halign{##\hskip0pt\cr\rm lim\cr
	\noalign{\vskip1pt}
	$\scriptstyle\mathord-\mskip-10mu plus1fil
	\mathord-\mskip-10mu plus1fil
	\mathord\rightarrow$\cr}}}\nolimits}
\def\lng{\mathop{\rm length}\nolimits}
\def\Max{\mathop{\rm Max}\nolimits}
\def\mod#1{\ifinner\mskip8mu(\mathop{\rm mod}#1)
        \else\mskip12mu(\mathop{\rm mod}#1)\fi}
\def\op{^{\mathop{\rm op}}}
\def\Spec{\mathop{\rm Spec}}
\def\lmapr#1#2{{}\mathrel{\smash{\mathop{\count0=#1
  \loop
    \ifnum\count0>0
    \advance\count0 by-1\smash{\mathord-}\mkern-4mu
  \repeat
  \mathord\rightarrow}\limits^{#2}}}{}}
\let\al\alpha
\let\ep\varepsilon
\let\la\lambda
\let\ph\varphi
\let\th\theta
\let\ze\zeta
\let\De\Delta
\def\0{_{(0)}}
\def\1{_{(1)}}
\def\2{_{(2)}}

\def\a{{\frak a}}
\def\at{\tilde a}
\def\b{{\frak b}}
\def\eh{\hat e}
\def\m{{\frak m}}
\def\mh{\hat m}
\def\mt{\tilde m}
\def\n{{\frak n}}
\def\p{{\frak p}}
\def\q{{\frak q}}
\def\vt{\tilde v}
\def\xt{\tilde x}
\def\C{{\cal C}}
\def\I{{\cal I}}
\def\F{{\cal F}}
\def\G{{\cal G}}
\def\H{{\cal H}}
\def\Hd{{H^\circ}}
\def\M{{\cal M}}

\def\U{{\cal U}}
\def\Ut{{\widetilde U}}
\def\Z{{\Bbb Z}}
\def\Q{{\Bbb Q}}
\def\AH{\HM_A}
\def\AM{\hbox{}_A\mskip1mu\M}
\def\RM{\hbox{}_R\mskip1mu\M}
\def\FAH{\F_A^H}
\def\MAH{\M_A^H}
\def\MBH{\M_B^H}
\def\HM{\hbox{}_H\M}
\def\HdM{\hbox{}_\Hd\M}
\def\HdMA{\HdM_A}
\def\KG{K\mkern1mu\backslash\mkern1mu G}
\def\ds{\le}
%
%		Start of the paper
%
\frontmatter

\title{Projectivity of Hopf algebras over subalgebras\break
with semilocal central localizations}
\author{Serge Skryabin}
\address{Chebotarev Research Institute,
Universitetskaya St.~17, 420008 Kazan, Russia}
\email{Serge.Skryabin@ksu.ru}

\endfrontmatter

\section
Introduction

Let $H$ be a Hopf algebra over a field $k$ and $A$ a {\it right coideal 
subalgebra} of $H$, that is, $A$ is a subalgebra satisfying $\De(A)\sbs A\ot 
H$ where $\De$ is the comultiplication in $H$. In case when $H$ is finitely 
generated commutative, the right coideal subalgebras are intimately related to 
the homogeneous spaces for the corresponding group scheme. The purpose of this 
paper is to extend the class of pairs $A$, $H$ for which $H$ is proved to be 
either projective or flat as a module over $A$. As is known the faithful 
flatness over Hopf subalgebras may be lacking in general. Examples given by 
Schauenburg \cite{Scha00} use some extremely big Hopf algebras coming from a 
universal construction of \cite{Tak71}. Positive results can be expected 
therefore only under some finiteness assumptions.

A Hopf algebra is called {\it residually finite dimensional} \cite{Mo} if its 
ideals of finite codimension have zero intersection. Many important classes of 
Hopf algebras are residually finite dimensional. Among them are the finitely 
generated commutative Hopf algebras, the universal enveloping algebras of 
finite dimensional Lie algebras, and also Hopf algebras related to quantum 
groups.

We say that a ring $R$ has {\it semilocal localizations with respect to 
a central subring $Z$} if for each maximal ideal $\m$ of $Z$ the localization 
$R_\m$ of $R$ at the multiplicatively closed set $Z\setm\m$ is a semilocal 
ring whose Jacobson radical contains $\m R_\m$. For instance, this property is 
satisfied for any ring module-finite over a central subring. For each ring $R$ 
let $\M_R$ and $\RM$ denote the categories of right and left $R$-modules, 
respectively.

\proclaim
Theorem 0.1.
Let $H$ be a residually finite dimensional Hopf algebra{\rm,} and let $A$ be a 
Hopf subalgebra having semilocal localizations with respect to a central 
subring $Z$. Then $H$ is a projective generator in $\M_A$ and in $\AM$.
\endproclaim

\proclaim
Theorem 0.2.
Let $A\sbs B\sbs H$ where $H$ is a residually finite dimensional Hopf algebra,
$B$ is a Hopf subalgebra, and $A$ is a right coideal subalgebra having
semilocal localizations with respect to a central subring $Z$. Suppose
$B$ is right module-finite over $A$ and the antipode of $B$ is bijective.
Then $H$ is a projective generator in $\M_A$.
\endproclaim

In both theorems we encounter projective modules of a very special kind. In 
fact $H\ot_AA_\m$ is a free $A_\m$-module for any maximal ideal $\m$ of $Z$. 
When $\dim H<\infty$, Theorem 0.2 applies to an arbitrary right coideal 
subalgebra $A$ since we may take $Z=k$. In this case $H$ is a free 
$A$-module, which generalizes the Nichols-Zoeller freeness theorem 
\cite{Ni89}. The investigation of the freeness over right coideal subalgebras 
in the finite dimensional case was initiated in \cite{Hoff}, \cite{Ma92}, and 
the full solution was obtained in \cite{Sk07}.

Over coideal subalgebras one can expect a Hopf algebra to be a flat module 
rather than faithfully flat or projective. If $H=k[G]$ is the function algebra 
of an affine group scheme $G$ of finite type over $k$ then for any group 
subscheme $K$ of $G$ the function algebra $A=k[\KG]$ on the right homogeneous 
space $\KG$ is a right coideal subalgebra of $H$. There are many cases where
$\KG$ is quasiaffine, and so $\KG$ may be identified with an open 
subscheme $U$ of the affine scheme $\Spec A$. Since the canonical morphism
$G\to\KG$ is flat, $H$ is always flat over $A$. However, the faithful 
flatness and projectivity are obtained precisely when $U=\Spec A$, i.e., when 
$\KG$ is affine.

\proclaim
Theorem 0.3.
Let $A\sbs B\sbs H$ where $H$ is a directed union of residually 
finite dimensional Hopf subalgebras, $B$ is any Hopf subalgebra, and $A$ is a 
right coideal subalgebra contained in the center of $B$.
Then $H$ is flat in $\M_A$. If $A$ is a Hopf subalgebra then $H$ is a 
projective generator in $\M_A$ and in $\AM$.
\endproclaim

Any commutative Hopf algebra $H$ is a directed union of finitely generated 
Hopf subalgebras, and those are residually finite dimensional. In this case
Theorem 0.3 applies to an arbitrary right coideal subalgebra $A$; we recover
the projectivity result of Takeuchi \cite{Tak79} and the flatness result of 
Masuoka and Wigner \cite{Ma94}. Our result is new even when $A$ is a central 
Hopf subalgebra of $H$. An interesting known example is the quantized function 
algebra at a root of unity; this Hopf algebra contains the ordinary function 
algebra of a semisimple algebraic group in its center. Projectivity was proved
in this case by De Concini and Lyubashenko \cite{Con94}; they needed detailed
information about quantized function algebras.

Several related results are known where $H$ is not assumed to be residually 
finite dimensional, but there are restrictions of a different kind. As was 
established by Schneider \cite{Schn93}, any left or right noetherian Hopf 
algebra is a faithfully flat module over central Hopf subalgebras. More 
recently Wu and Zhang \cite{Wu03} discovered that the projectivity holds for
finite extensions of finitely generated PI Hopf algebras under certain 
finiteness assumptions about injective or projective dimensions. Of a somewhat
different flavor are results for pointed Hopf algebras \cite{Kh}, \cite{Rad77}
or Hopf algebras with cocommutative coradical \cite{Ma91} which impose a 
restriction on the coalgebra structure rather than the algebra structure.

As was emphasized in \cite{Sk07}, it is natural to investigate projectivity
in the more general settings where $A$ is assumed to be a (right) 
{\it$H$-comodule algebra}. Such an algebra $A$ has a right $H$-comodule 
structure given by an algebra homomorphism $\rho_A:A\to A\ot H$. With $A$ one 
associates the category of {\it right Hopf modules} $\MAH$ \cite{Doi83},
\cite{Tak79}. The objects of $\MAH$ have structures of a right $A$-module and 
a right $H$-comodule such that the comodule structure map $\rho_M:M\to M\ot H$ 
becomes an $\M_A$-morphism if we let $A$ operate on $M\ot H$ via $\rho_A$. 
When $A$ and $H$ are commutative, $\MAH$ is equivalent to the category of 
$G$-linearized quasicoherent sheaves on $\Spec A$, where $G$ is the group 
scheme corresponding to $H$. As usual we use the term ``{\it$H$-costable}" in 
the sense ``stable under the coaction of $H$".

\proclaim
Problem 1.
Let $A$ be an {\rm$H$-simple $H$-comodule algebra}, i.e., $A$ has no
$H$-costable ideals other than $0$ and $A$. For what classes of algebras is 
every nonzero object of $\MAH$ a projective generator in $\M_A${\rm?} 
\endproclaim

When $H=k$ is the trivial Hopf algebra, the $H$-simplicity of $A$ means that 
$A$ is simple, and the question asks whether all right $A$-modules are
projective. Certainly, this holds if and only if $A$ is artinian. In general
Problem 1 is meaningful under the assumption that $A$ has a simple artinian 
factor ring. This is automatic for subalgebras of $H$ since the counit 
$\ep:H\to k$ makes $k$ into a factor algebra of any subalgebra. In order to 
treat flatness we want to weaken the assumption about the $H$-simplicity of 
$A$.

\proclaim
Problem 2.
Let $A$ be an $H$-costable subalgebra of an $H$-comodule algebra $B$. Suppose 
that $IB=B$ for each nonzero $H$-costable ideal $I$ of $A$. For what classes 
of algebras is every object of $\MBH$ flat in $\M_A${\rm?}
\endproclaim

There is a dual formulation for (left) {\it$H$-module algebras}. Here $A$ is 
an algebra which has a left $H$-module structure compatible with the 
multiplication; $\MAH$ is replaced with the category $\HM_A$ whose objects 
have a right $A$-module structure and a compatible left $H$-module structure.  
Working with module algebras gives some advantage since in this case the 
coalgebra structure on $H$ is important, and we can use the family of finite 
dimensional subcoalgebras. In fact we are only able to approach Problems 1 and 
2 for $H$-comodule algebras by making a reduction to similar questions for 
module algebras over the finite dual $\Hd$ of $H$. The correct correspondence
between the $H$-comodule structure and the $\Hd$-module structure is available
when $H$ is residually finite dimensional. This explains why this kind of 
restriction on $H$ appears in Theorems 0.1--0.3.

An object $M\in\HM_A$ is called {\it$A$-finite} if $M$ is finitely generated 
as an $A$-module; $M$ is {\it locally $A$-finite} if $M$ is a directed union 
of $A$-finite subobjects. Theorem 5.6 provides a projectivity result for 
locally $A$-finite objects of $\HM_A$ assuming that $A$ has semilocal 
localizations with respect to a central subring. This unifies the cases of 
commutative algebras and semilocal ones considered in \cite{Sk07}.

The short proof of the previous result proposed in \cite{Sk07} for the case
of commutative $A$ and cocommutative $H$ is based on three properties of the 
Fitting ideals of a finitely generated $A$-module $M$: (1) these ideals 
contain enough information to recognize projective modules of constant rank,
(2) they behave functorially with respect to the change of ring, (3) they are 
stable under a compatible action of $H$. The definition of the Fitting ideals
involves computing determinants, so it does not generalize to noncommutative 
rings.

In section 1 of the present article we introduce certain ideals $I_r(M)$ of a 
ring $R$ for any rational number $r\ge0$ and a finitely generated right 
$R$-module $M$ imposing some assumptions about the localizations $R_\m$ at the 
maximal ideals of a central subring of $R$. There are analogs, though less 
satisfactory, of the three properties mentioned above. When $R=A$ is an 
$H$-module algebra with semilocal central localizations and $M\in\HM_A$, the 
ideals $I_r(M)$ are not $H$-stable in general. Our expectation is that the 
normalized rank $r_P(M)$ at a maximal ideal $P$ of $A$ is determined by those 
rational values $r$ for which $P$ contains the smallest $H$-stable ideal 
$J_r(M)$ of $A$ such that $I_r(M)\sbs J_r(M)$. When this holds, 
$r_P(M)=r_Q(M)$ for any pair $P,Q$ of maximal ideals of $A$ containing the 
same $H$-stable ideals. This property can be viewed as the $H$-invariance of 
the rank function $P\mapsto r_P(M)$ defined on the maximal spectrum $\Max A$ 
of $A$. We are able to prove it only under some technical restrictions. As a 
result, we gain less control over the situation in those cases where $A$ is 
not $H$-simple. This leads to more restrictive assumptions about $A$ when 
dealing with flatness.

We will use standard notation from the theory of Hopf algebras \cite{Mo}, 
\cite{Sw}. For each ring $R$ denote by $\Jac(R)$ the Jacobson radical of $R$, 
by $\Max R$ and $\Spec R$ the maximal and prime spectra of $R$. A ring $R$ is
{\it semilocal} if $R/\Jac(R)$ is artinian; a semilocal ring with a single 
maximal ideal is {\it quasilocal}. A semilocal ring $R$ is {\it semiprimary}
if $\Jac(R)$ is nilpotent; a semiprimary ring with a single maximal ideal is 
{\it primary}. Denote by $\Z_+$ the semigroup of positive integers. 

\section
1. Construction of ideals

Let $R$ be a ring. If $M\in\M_R$ is generated by elements $e_1,\ldots,e_n$, we 
denote by $\I_{e_1,\ldots,e_n}$ the ideal of $R$ generated by all elements of 
$R$ which occur as a coefficient in a zero linear combination 
$e_1x_1+\cdots+e_nx_n=0$ with $x_1,\ldots,x_n\in R$.

\proclaim
Lemma 1.1.
Suppose that $e_1,\ldots,e_n$ generate $M\in\M_R$. If $\ph:R\to R'$ is a ring 
homomorphism and $e'_i=e_i\ot1\in M\ot_RR'$, then 
$\I_{e'_1,\ldots,e'_n}=R'\ph(\I_{e_1,\ldots,e_n})R'$.
\endproclaim

\Proof.
Let $\th:R^n\to M$ be the epimorphism in $\M_R$ sending the standard 
generators of the free module $R^n$ to $e_1,\ldots,e_n$. For $i=1,\ldots,n$ 
denote by $\pi_i:R^n\to R$ the projection onto the $i$th summand. 
The ideal $\I_{e_1,\ldots,e_n}$ is generated by $\sum\pi_i(K)$ where 
$K=\Ker\th$. Tensoring with $R'$, we get an exact sequence of $R'$-modules
$$ 
K\ot_RR'\lmapr3{\al\ot\id}R^n\ot_RR'\lmapr3{\th\ot\id}M\ot_RR'\to0
$$ where 
$\al:K\to R^n$ is the inclusion map. We have an isomorphism 
$R^n\ot_RR'\cong R'^n$ with $\pi_i\ot\id:R^n\ot_RR'\to R'$ giving the 
projection onto the $i$th summand. Hence $\I_{e'_1,\ldots,e'_n}$ 
coincides with the ideal of $R'$ generated by $\sum\,(\pi_i\ot\id)(K')$ where 
$K'=\Ker(\th\ot\id)$. Since $K'$ is equal to the image of $\al\ot\id$, the 
conclusion is clear.
\endproof

Recall that $R$ is said to be {\it weakly finite} if for each integer $n>0$
every generating set for the free right $R$-module $R^n$ containing exactly 
$n$ elements is a basis for $R^n$; equivalently, every $R$-module epimorphism
$R^n\to R^n$ is an isomorphism. This can also be reformulated in terms of
invertibility of $n\times n$-matrices with entries in $R$.

\proclaim
Lemma 1.2.
Suppose that all factor rings of $R$ are weakly finite. If $e_1,\ldots,e_n$ 
and $e'_1,\ldots,e'_n$ are two systems of generators for $M$ having the same 
number of elements then $\I_{e'_1,\ldots,e'_n}=\I_{e_1,\ldots,e_n}$.  
\endproclaim

\Proof.
This follows from \cite{Sk07, Lemma 2.3}. By that lemma
$\I_{e'_1,\ldots,e'_n}\sbs\I_{e_1,\ldots,e_n}$ since 
$R/\I_{e_1,\ldots,e_n}$ is weakly $n$-finite. The opposite inclusion holds by 
symmetry.
\endproof

If $R$ has weakly finite factor rings and $M$ is $n$-generated then we put
$$
I_n(M)=\I_{e_1,\ldots,e_n}
$$
where $e_1,\ldots,e_n$ is any set of $n$ generators for $M$. By Lemma 1.2
the above ideal does not depend on the choice of a generating set. The ideals
$I_r(M)$ are thus defined for all integers $r\ge\mu(M)$ where $\mu(M)$ denotes
the minimal number of generators for $M$. When $M=0$ we put $\mu(M)=0$ and
$I_0(M)=0$ for consistency reasons.

We do not indicate the base ring explicitly in the notation for $I_r(M)$. 
Given a ring homomorphism $R\to R'$, let $I_r(M\ot_RR')$ be the ideal of $R'$ 
corresponding to the induced $R'$-module $M\ot_RR'$ (when defined). 
Especially, this convention will be in force when $R'$ is either a factor ring 
or an Ore localization of $R$.

\setitemsize(iii)
\proclaim
Lemma 1.3.
Suppose that $R$, $R'$ are two rings with weakly finite factor rings,
$\ph:R\to R'$ a homomorphism and $M$, $N$ two finitely generated right 
$R$-modules.

\item(i)
$I_r(M)=R$ for all integers $r>\mu(M)$.

\item(ii)
$M\cong R^n$ if and only if $M$ is $n$-generated with $I_n(M)=0$.

\item(iii)
$I_{r+s}(M\oplus N)=I_r(M)+I_s(N)$ for all integers $r\ge\mu(M)$ and 
$s\ge\mu(N)$.

\item(iv)
$I_{rt}(M^t)=I_r(M)$ for all integers $r\ge\mu(M)$ and $t>0$.

\item(v)
$I_r(M\ot_RR')=R'\ph\bigl(I_r(M)\bigr)R'$ for all integers $r\ge\mu(M)$.

\endproclaim

\Proof.
(i) If $r>\mu(M)$ then $M$ can be generated by $r-1$ elements, say 
$e_1,\ldots,e_{r-1}$. Adding another element $e_r=0$, we get a set of $r$ 
generators for $M$. Now $I_r(M)=\I_{e_1,\ldots,e_r}$.  However, 
$1\in\I_{e_1,\ldots,e_r}$ since $1$ is a coefficient in the relation $e_r=0$.

(ii) Let $M$ be generated by $e_1,\ldots,e_n$. Clearly, $e_1,\ldots,e_n$ is a 
basis for $M$ if and only if $\I_{e_1,\ldots,e_n}=0$.

(iii) Let us identify $M$ and $N$ with submodules of $M\oplus N$. Pick 
generating sets $e_1,\ldots,e_r$ for $M$ and $e'_1,\ldots,e'_s$ for $N$. Then 
the $r+s$ elements $e_1,\ldots,e_r,e'_1,\ldots,e'_s$ generate $M\oplus N$. 
Given $x_1,\ldots,x_r,y_1,\ldots,y_s\in R$, the equality
$\sum e_ix_i+\sum e'_jy_j=0$ holds if and only if both $\sum e_ix_i=0$ and 
$\sum e'_jy_j=0$. Hence
$$
\I_{e_1,\ldots,e_r,e'_1,\ldots,e'_s}=\I_{e_1,\ldots,e_r}+\I_{e'_1,\ldots,e'_s}.
$$

(iv) This follows from (iii) by induction on $t$.

(v) This is a restatement of Lemma 1.1.
\endproof

Suppose that $M\in\M_R$ is finitely generated and $r\ge0$ is a rational number.
The set $\{l\in\Z_+\mid rl\ge\mu(M^l)\}$
is closed under addition since $\mu(M^{l+t})\le\mu(M^l)+\mu(M^t)$ for all 
$l,t\in\Z_+$. If the inequality $rl\ge\mu(M^l)$ holds for at least one $l$ 
then we can find such an $l$ with the property that $rl\in\Z$, replacing $l$ 
with a suitable multiple if necessary. We put then
$$
I_r(M)=I_{rl}(M^l),
$$
which does not depend on the choice of an $l$ with the above properties. In 
fact, if $t\in\Z_+$ also satisfies $rt\in\Z$ and $rt\ge\mu(M^t)$ then
$I_{rl}(M^l)=I_{rlt}(M^{lt})=I_{rt}(M^t)$ by (iv) of Lemma 1.3. If $r\in\Z$
and $r\ge\mu(M)$ then the initial definition of $I_r(M)$ agrees with the newer 
one since $l=1$ satisfies the required properties.

We say that $r$ is {\it$M$-admissible} if $rl\ge\mu(M^l)$ for some
$l\in\Z_+$. We have defined the ideals $I_r(M)$ for all $M$-admissible 
rational numbers. All statements in the next lemma immediately reduce to the 
corresponding statements in Lemma 1.3.

\proclaim
Lemma 1.4.
Retaining the assumptions about $R$, $R'$, $\ph$, $M$, $N$ as in Lemma 
{\rm1.3,} let $r\in\Q$ be $M$-admissible and $s\in\Q$ be $N$-admissible.

\item(i)
$I_r(M)=R$ whenever $rl>\mu(M^l)$ for some $l\in\Z_+$.

\item(ii)
$M^l\cong R^n$ for integers $l>0$, $n\ge0$ if and only if\/ $n\ge\mu(M^l)$ and
$I_{n/l}(M)=0$.

\item(iii)
$I_{r+s}(M\oplus N)=I_r(M)+I_s(N)$.

\item(iv)
$I_{rt}(M^t)=I_r(M)$ for all $t\in\Z_+$.

\item(v)
$I_r(M\ot_RR')=R'\ph\bigl(I_r(M)\bigr)R'$.

\endproclaim

Part (v) of this lemma is valid in a slightly more general situation where $R$ 
is not assumed to have weakly finite factor rings.

\proclaim
Lemma 1.5.
Suppose that $r\in\Q$ is $M$-admissible. Then $R$ has an ideal $K$ such that 
$I_r(M\ot_RR')=R'\ph(K)R'$ for any homomorphism $\ph:R\to R'$ into a ring $R'$ 
with weakly finite factor rings.
\endproclaim

\Proof.
There exists $l\in\Z_+$ such that $n=rl\in\Z$ and $n\ge\mu(M^l)$. Then $M^l$ 
is generated by $n$ elements, say $e_1,\ldots,e_n$. Take 
$K=\I_{e_1,\ldots,e_n}$. For any $\ph$ satisfying the hypotheses the 
$R'$-module $M^l\ot_RR'$ is generated by $e'_1,\ldots,e'_n$ where 
$e'_i=e_i\ot1$ for each $i$. Lemma 1.1 yields
$I_r(M\ot_RR')=I_n(M^l\ot_RR')=\I_{e'_1,\ldots,e'_n}=R'\ph(K)R'$.  
\endproof

If $M\in\M_R$ is finitely generated and $P\in\Max R$ is such that $R/P$ is 
simple artinian, then we put
$$
r_P(M)={\lng M/MP\over\lng R/P}
$$
where $\lng$ stands for the composition series length in $\M_R$.

If $R$ is semilocal then the set $\Max R$ is finite and coincides with the set 
of primitive ideals of $R$. For each $P\in\Max R$ the ring $R/P$ is simple 
artinian, so that $r_P(M)$ is defined. Recall that any semilocal ring is 
weakly finite. Moreover, all factor rings of such a ring are themselves 
semilocal, hence weakly finite.

\setitemsize(ii)
\proclaim
Lemma 1.6.
Suppose $R$ is semilocal, $M\in\M_R$ is finitely generated and $r\in\Q$.

\item(i)
$\mu(M)=\min\{n\in\Z\mid n\ge r_P(M)\hbox{ for all }P\in\Max R\}$.

\item(ii)
$r$ is $M$-admissible if and only if $r\ge r_P(M)$ for all $P\in\Max R$.

\endproclaim

\Proof.
(i) Let $J=\Jac(R)$. We have $M/M\!J\cong\prod_{P\in\Max R}M/MP$ since $R/J$ 
is semisimple artinian. It follows from Nakayama's Lemma that $M$ is 
$n$-generated if and only if so is $M/MP$ for each $P$. Since every 
$R/P$-module is isomorphic to a direct sum of copies of the simple module, 
$M/MP$ is an epimorphic image of $(R/P)^n$ if and only if $\,\lng 
M/MP\le\lng\,(R/P)^n$, which can be rewritten as $r_P(M)\le n$. 

(ii) Put $m=\max\{r_P(M)\mid P\in\Max R\}$. It was proved in (i) that 
$\mu(M^l)$ is equal to the smallest integer $n$ such that
$n\ge r_P(M^l)=r_P(M)l$ for all $P$, i.e., $n\ge ml$. Thus $rl\ge\mu(M^l)$
implies $r\ge m$. If $r=m$ then the required inequality holds for any 
$l\in\Z_+$ such that $ml\in\Z$.
\endproof

Let $Z\sbs R$ be a central subring. Denote by $R_z$ and $R_\p$, respectively, 
the localizations of $R$ at the multiplicatively closed subsets 
$\{z^i\mid i=0,1,\ldots\}$ and $Z\setm\p$ where $z\in Z$ is any element 
and $\p$ a prime ideal of $Z$. Similarly, $M_z$ and $M_\p$ will denote the 
respective localizations of $M\in\M_R$. If $M$ is finitely generated, then 
$M_z$ is a finitely generated $R_z$-module and $M_\p$ a finitely generated 
$R_\p$-module. For each rational number $r\ge0$ put
$$
T_r(M)=\{z\in Z\mid r\hbox{ is $M_z$-admissible}\},
$$\removelastskip
$$
\Ut_r(M)=\{\p\in\Spec Z\mid r\hbox{ is $M_\p$-admissible}\},
$$\removelastskip
$$
U_r(M)=\{\m\in\Max Z\mid r\hbox{ is $M_\m$-admissible}\}.
$$
The open subsets $D(z)=\{\p\in\Spec Z\mid z\notin\p\}$ with $z\in Z$ give 
a basis for the topology on $\Spec Z$. If $\p\in D(z)$ then $R_\p$ is a 
localization of $R_z$ and $M_\p\cong M_z\ot_{R_z}R_\p$, whence 
$\mu(M_\p^l)\le\mu(M_z^l)$ for each $l\in\Z_+$. It follows that 
$D(z)\sbs\Ut_r(M)$ whenever $z\in T_r(M)$.

\proclaim
Lemma 1.7.
Any $\p\in\Ut_r(M)$ is contained in $D(z)$ for some $z\in T_r(M)$. Hence 
$\Ut_r(M)$ is open in $\Spec Z$ and $U_r(M)=\Ut_r(M)\cap\Max Z$ is open in 
$\Max Z$.
\endproclaim

\Proof.
Let $l>0$ be an integer such that $n=rl\in\Z$ and $n\ge\mu(M_\p^l)$. We can 
find $n$ elements $e_1,\ldots,e_n\in M^l$ whose images in $M_\p^l$ 
generate the latter $R_\p$-module. For each $x\in M^l$ there exists $z\in 
Z\setm\p$ such that $xz$ is an $R$-linear combination of 
$e_1,\ldots,e_n$. As $M^l$ is finitely generated, we can find a $z$ which 
fulfills the required property for all $x$ simultaneously. Then the 
$R_z$-module $M_z^l$ is generated by the images of $e_1,\ldots,e_n$, whence 
$n\ge\mu(M_z^l)$. Thus $z\in T_r(M)$ and $D(z)$ is a neighborhood of $\p$ 
contained in $\Ut_r(M)$.
\endproof

\proclaim
Lemma 1.8.
Suppose that $\m R_\m\sbs\Jac(R_\m)$ for each $\m\in\Max Z$. Then any finite 
subset $X\sbs U_r(M)$ is contained in $D(z)$ for some $z\in T_r(M)$.  
\endproclaim

\Proof.
We can find $l\in\Z_+$ such that $n=rl\in\Z$ and $n\ge\mu(M_\m^l)$ for all 
$\m\in X$. Since $M_\m/M_\m\m\cong M/M\m$, we have $n\ge\mu(M^l\!/M^l\m)$ for 
$\m\in X$. Since $\m+\n=Z$ for any pair of distinct ideals $\m,\n\in X$, the
canonical map $M^l\to\prod_{\m\in X}M^l/M^l\m$ is surjective by Chinese 
Remainder Theorem. There exist $e_1,\ldots,e_n\in M^l$ whose cosets modulo 
$\m$ generate the $R$-module $M^l\!/M^l\m$ for each $\m\in X$. Then for each 
$\m\in X$ the images of $e_1,\ldots,e_n$ in $M_\m^l$ generate the 
$R_\m$-module $M_\m^l$ modulo $M_\m^l\m$; by Nakayama's Lemma $M_\m^l$ is 
generated by those images.

For $x\in M^l$ denote by $\a_x$ the ideal of $Z$ consisting of those elements 
$z\in Z$ for which $xz$ lies in the submodule $N$ of $M^l$ generated by 
$e_1,\ldots,e_n$. Since $N_\m=M^l_\m$, we have $\a_x\not\sbs\m$ for any
$\m\in X$; hence $\a_x\not\sbs\bigcup_{\m\in X}\m\,$ \cite{Bou, Ch. II, \S1, 
Prop. 2}. So there exists $z\in Z$ which lies in none of the ideals from $X$ 
and satisfies $xz\in N$. As $M^l$ is finitely generated, we can find a $z$ 
which fulfills that property for all $x$ simultaneously. We obtain $z\in 
T_r(M)$ and $X\sbs D(z)$.
\endproof

Further on in this section we make the following assumption:

\medskip
\noindent(A)
{\it
$R_\m$ has weakly finite factor rings and $\m R_\m\sbs\Jac(R_\m)$ for each
$\m\in\Max Z$.
}

\medskip
Since $R$ is embedded in $\prod_{\m\in\Max Z}R_\m$, the weak finiteness of $R$
follows from the weak finiteness of all localizations $R_\m$. This 
observation, applied to the factor rings of $R$, shows that $R$ has weakly 
finite factor rings provided that so do all rings $R_\m$.

Let $M\in\M_R$ be finitely generated and $r\in\Q$. When $r$ is $M$-admissible,
we have $U_r(M)=\Max Z$ and $I_r(M_\m)=I_r(M)R_\m$ for all $\m\in\Max Z$ by 
Lemma 1.4(v). Every ideal of $R$ is completely determined by its extensions to 
the rings $R_\m$. In particular, $I_r(M)$ consists precisely of those elements 
$a\in R$ whose image $a_\m$ in $R_\m$ belongs to $I_r(M_\m)$ for each 
$\m\in\Max Z$.

We now extend the range of $r$ in the definition of $I_r(M)$ to arbitrary 
nonnegative values. The ideal $I_r(M_\m)$ of $R_\m$ has already been defined 
when $\m\in U_r(M)$. Put
$$
I_r(M)=\{a\in R\mid a_\m\in I_r(M_\m){\rm\ for\ each\ }\m\in U_r(M)\}.
$$
Clearly $I_r(M)$ is an ideal of $R$. If $U_r(M)=\varnothing$ then 
$I_r(M)=R$.

\proclaim
Lemma 1.9.
If $U_r(M)$ is quasicompact, then $I_r(M)R_\m=I_r(M_\m)$ for $\m\in U_r(M)$.
\endproclaim

\Proof.
The ideals of $R_\m$ are extensions of ideals of $R$. Fixing $\m$, we have to 
show that for each $a\in R$ with $a_\m\in I_r(M_\m)$ there exists $s\in 
Z\setm\m$ such that $as\in I_r(M)$, i.e., $a_\n s_\n\in I_r(M_\n)$ for all
$\n\in U_r(M)$.  

Suppose that $\m\in D(z)$ where $z\in T_r(M)$. As $r$ is $M_z$-admissible,
Lemma 1.5 can be applied to the $R_z$-module $M_z$. Let $K$ be the ideal of 
$R_z$ given by that lemma. For each $\n\in D(z)\cap\Max Z$ the ring $R_\n$ is 
a localization of $R_z$, whence $I_r(M_\n)=KR_\n$ by Lemma 1.5. Applying this 
formula with $\m=\n$, we deduce that $a_zt_z\in K$ for a suitable $t\in 
Z\setm\m$ where $a_z,t_z$ denote the images of $a,t$ in $R_z$. Passing now to 
$R_\n$, we see that $a_\n t_\n\in I_r(M_\n)$ for any $\n$ as above.

Given an arbitrary $\n\in U_r(M)$, there exists $z\in T_r(M)$ such that 
$\{\m,\n\}\sbs D(z)$ by Lemma 1.8. Hence $U_r(M)$ is covered by the open 
subsets $D(z)$ with $z\in T_r(M)$, $\,z\notin\m$. Since $U_r(M)$ is 
quasicompact, we have $U_r(M)\sbs D(z_1)\cup\cdots\cup D(z_n)$ for some
elements $z_1,...,z_n\in T_r(M)$ such that $\m\in D(z_i)$ for each 
$i=1,\ldots,n$. We have seen that for each $i$ there exists $t_i\in Z\setm\m$ 
such that $a_\n(t_i)_\n\in I_r(M_\n)$ for all $\n\in D(z_i)\cap\Max Z$. Now
$s=t_1\cdots t_n$ is the desired element.
\endproof

\Remark.
The topological space $\Max Z$ is always quasicompact. Hence Lemma 1.9 applies 
for any $r$ with $U_r(M)=\Max Z$. If $\Max Z$ is noetherian (e.g., if 
$Z/\Jac(Z)$ is noetherian), then every open subset of $\Max Z$ is 
quasicompact. In this case any nonnegative value of $r$ is legitimate.  
\endremark

\proclaim
Lemma 1.10.
Let $\ph:R\to R'$ be a ring homomorphism where $R$ satisfies {\rm(A)}, while 
$R'$ has weakly finite factor rings. Suppose that $U_r(M)$ is quasicompact and
there exists a finite subset $X\sbs U_r(M)$ such that $\ph(z)$ is invertible 
in $R'$ for each $z\in Z\setm\bigcup_{\m\in X}\m$. Then 
$I_r(M\ot_RR')=R'\ph\bigl(I_r(M)\bigr)R'$.
\endproclaim

\Proof.
By Lemma 1.8 there exists $z\in T_r(M)$ such that $X\sbs D(z)$, so that $z$ 
lies in none of the ideals $\m\in X$. Since $\ph(z)$ is invertible, $\ph$ 
extends to a homomorphism $\psi:R_z\to R'$. Recall that $r$ is 
$M_z$-admissible by the definition of $T_r(M)$. Let $K$ be the ideal of $R_z$ 
given by Lemma 1.5, when applied to the $R_z$-module $M_z$. Then 
$I_r(M\ot_RR')=R'\psi(K)R'$, and also $I_r(M_\m)=KR_\m$ for each $\m\in X$ as 
in the proof of Lemma 1.9. Thus the two ideals $K$ and $I_r(M)R_z$ of $R_z$ 
have the same extension to each ring $R_\m$ with $\m\in X$. Given any $a\in 
K$, there exists therefore $s\in Z\setm\bigcup_{\m\in X}\m$ such that $as_z\in 
I_r(M)R_z$ where $s_z$ denotes the image of $s$ in $R_z$. Since 
$\psi(s_z)=\ph(s)$ is invertible in $R'$, we deduce that 
$\psi(a)\in\psi\bigl(I_r(M)R_z\bigr)R'=\ph\bigl(I_r(M)\bigr)R'$.  Similarly, 
given any $a\in I_r(M)$, we prove that $\ph(a)\in\psi(K)R'$. Hence 
$R'\psi(K)R'=R'\ph\bigl(I_r(M)\bigr)R'$, and we are done.
\endproof

\Remark.
Lemma 1.10 will be used in the special case where $\ph$ is the canonical 
homomorphism onto a factor ring $R'$ of $R$. Suppose that $R'$ has finitely 
many maximal ideals and $\ph^{-1}(P)\cap Z\in U_r(M)$ for each 
$P\in\Max R'$. Take
$$X=\{\ph^{-1}(P)\cap Z\mid P\in\Max R'\}.
$$
If $z\in Z$, then $\ph(z)R'$ is an ideal of $R'$; hence either $\ph(z)R'=R'$ 
or $\ph(z)R'$ is contained in some $P\in\Max R'$. In the former case $\ph(z)$ 
is invertible, while in the latter case $z\in\m$ for some $\m\in X$. When 
$U_r(M)$ is quasicompact and $R'$ has weakly finite factor rings, 
the hypotheses of Lemma 1.10 are satisfied.
\endremark

Note that $U_r(M)\sbs U_s(M)$ whenever $r,s\in\Q$ satisfy $0\le r\le s$. Put
$$
\la(M)=\inf\{r\in\Q\mid U_r(M)=\Max Z\}.
$$

\setitemsize(iii)
\proclaim
Lemma 1.11.
Suppose $R$ satisfies {\rm(A)}, $M$ is finitely generated and $r\in\Q$,
$\,r\ge0$.

\item(i)
$I_r(M)=R$ whenever $r>\la(M)$.

\item(ii)
If $U_r(M)$ is quasicompact and $I_r(M)=0$ then for each $\m\in U_r(M)$ there 
exist integers $l>0$ and $n\ge0$ such that $r=n/l$ and $M^l_\m\cong R^n_\m$ in 
$\M_{R_\m}$.

\item(iii)
$I_{rt}(M^t)=I_r(M)$ for all $t\in\Z_+$.

\endproclaim

\Proof.
(i) Pick any $s\in\Q$ such that $\la(M)<s<r$. Then $U_s(M)=\Max Z$. For each 
$\m\in\Max Z$ there exists $l\in\Z_+$ such that $rl>sl\ge\mu(M_\m^l)$;
Lemma 1.4(i) shows that $I_r(M_\m)=R_\m$. The conclusion is now immediate from 
the definition of $I_r(M)$.

(ii) If $\m\in U_r(M)$, then there exists $l\in\Z_+$ such that $n=rl\in\Z$ and
$n\ge\mu(M^l_\m)$. By Lemma 1.9 $I_r(M_\m)=0$. Now we may apply Lemma 
1.4(ii).

(iii) Note that $U_{rt}(M^t)=U_r(M)$ by a straightforward check and 
$I_{rt}(M_\m^t)=I_r(M_\m)$ for each $\m\in U_r(M)$ by Lemma 1.4(iv).
\endproof

\Remark.
If the ring $R_\m$ is semilocal, then the isomorphism $M^l_\m\cong R^n_\m$ in 
(ii) holds for any pair of integers $l>0$ and $n\ge0$ such that $r=n/l$.
\endremark

\section
2. Rings with semilocal central localizations

We will assume throughout the whole section that $R$ has semilocal 
localizations with respect to a central subring $Z$. Thus
\medskip
\centerline{\it$R_\m$ is semilocal and\/ $\m R_\m\sbs\Jac(R_\m)$ \/for each\/
$\m\in\Max Z$.}
\medskip
\noindent
Note that $R$ satisfies assumption (A) from section 1. Any factor ring $R'$ of 
$R$ has semilocal localizations with respect to the image of $Z$ in $R'$.
In this section several properties of the ring $R$ will be stated for future 
use. Some of those are more or less known.

\proclaim
Lemma 2.1.
For any right primitive ideal $P$ of $R$ the ring $R/P$ is simple 
artinian and $P\cap Z\in\Max Z$. Given $\m\in\Max Z$, there are finitely many 
elements in the set
$$
\Max_\m R=\{P\in\Max R\mid P\cap Z=\m\}.
$$
The maximal ideals of $R_\m$ are precisely the ideals $P_\m=R_\m P$ with
$P\in\Max_\m R$.
\endproclaim

\Proof.
Let $P$ be the annihilator of a simple right $R$-module $V$. Since $1\notin P$,
there exists $\m\in\Max Z$ such that $P\cap Z\sbs\m$. The transformation 
$z_V$ of $V$ afforded by an element $z\in Z$ is an $\M_R$-endomorphism. Hence
the image and the kernel of $z_V$ are submodules of $V$. If $z\notin\m$, then 
$Vz\ne0$, whence $Vz=V$ and $\Ker z_V=0$ by the simplicity of $V$. In 
other words, $z_V$ is invertible for any $z\in Z\setm\m$. We may now regard 
$V$ as a simple $R_\m$-module. The condition $\m R_\m\sbs\Jac(R_\m)$ entails 
$V\m=0$, i.e., $\m\sbs P$. The maximality of $\m$ yields $P\cap Z=\m$. Now 
$P/\m R$ is a right primitive ideal of the factor ring $R/\m R\cong R_\m/\m 
R_\m$. The latter is semilocal since so is $R_\m$. Hence $R/\m R$ has 
finitely many primitive ideals, and the factor algebra by any of those is 
simple artinian.

If $P$ is any ideal of $R$ such that $\m\sbs P$ for some $\m\in\Max Z$, then 
$R_\m/P_\m\cong(R/P)_\m\cong R/P$; in this case $P_\m\in\Max R_\m$ if and only 
if $P\in\Max_\m R$. Any ideal $P'$ of $R_\m$ coincides with $P_\m$ where $P$ 
is the preimage of $P'$ in $R$; if $P'\in\Max R_\m$, then $\Jac(R_\m)\sbs P'$, 
and the assumption about $R_\m$ yields $\m\sbs P$, so that $P\in\Max_\m R$.
\endproof

Suppose further that $M$ is a finitely generated right $R$-module.

\proclaim
Lemma 2.2.
If $M\a=M$ for some ideal $\a$ of $Z$ then there exists $a\in\a$ such that 
$M(1-a)=0$. In particular, $\a=Z$ whenever $M$ is faithful.
\endproclaim

\Proof.
For each $\m\in\Max Z$ the right $R_\m$-module $M_\m$ is finitely generated and
$M_\m=M_\m\a$. If $\a\sbs\m$, then $\a R_\m\sbs\Jac(R_\m)$, which forces 
$M_\m=0$ by Nakayama's Lemma. In this case any element of $M$ is annihilated by 
some element in $Z\setm\m$; since $M$ is finitely generated over $R$ and $Z$ 
is in the center of $R$ there exists $z\in Z\setm\m$ such that $Mz=0$. Denote 
by $\b$ the annihilator of $M$ in $Z$. We conclude that $\a+\b$ cannot be 
contained in any maximal ideal of $Z$, whence $\a+\b=Z$. It follows that 
$1-a\in\b$ for some $a\in\a$. If $M$ is faithful, we must have $a=1$.  
\endproof

\proclaim
Lemma 2.3.
For any $r\in\Q$, $\,r\ge0$, we have
$$
U_r(M)=\{\m\in\Max Z\mid r_P(M)\le r\hbox{ for all }P\in\Max_\m R\}.
$$
\endproclaim

\Proof.
Let $\m\in\Max Z$. Lemma 1.6(ii) together with Lemma 2.1 and the definition of 
$U_r(M)$ in section 1 show that $\m\in U_r(M)$ if and only if 
$r_{P_\m}(M_\m)\le r$ for all $P\in\Max_\m R$. Since the image of $Z\setm\m$ 
in $R/P$ consists of invertible elements, we have $R_\m/P_\m\cong R/P$ and 
$M_\m/M_\m P_\m\cong M/MP$, whence $r_{P_\m}(M_\m)=r_P(M)$.
\endproof

\proclaim
Lemma 2.4.
The supremum $r(M)=\sup\{r_P(M)\mid P\in\Max R\}$
is attained at some maximal ideal of $R$.
\endproclaim

\Proof.
Note that $\Max Z\ne U_s(M)$ for any $s\in\Q$ such that $0\le s<r(M)$. 
Indeed, if $P\in\Max R$ satisfies $r_P(M)>s$, then $P\cap Z$ is a maximal 
ideal of $Z$ lying outside of $U_s(M)$ by Lemma 2.3. Since $\Max Z$ is 
quasicompact and each $U_s(M)$ is open in $\Max Z$, we get
$$
\Max Z\ne\bigcup_{0\le s<r(M)}U_s(M).
$$
Pick $\m\in\Max Z$ contained in none of the subsets $U_s(M)$ with $s<r(M)$.  
Next, in the finite set $\Max_\m R$ pick $P$ with the maximum value of 
$r_P(M)$. By Lemma 2.3 $r_P(M)>s$ for any $s\in\Q$ with $s<r(M)$. Hence 
$r_P(M)=r(M)$.
\endproof

\proclaim
Lemma 2.5.
Let $K$ be an ideal of $R$ such that $R/K$ is semilocal, and let $r=n/l$ 
for some integers $n\ge0$, $l>0$. Suppose that $U_r(M)$ is quasicompact and 
$r_Q(M)\le r$ for each $Q\in\Max R$ such that there exists $P\in\Max R$ 
satisfying $P\sps K$ and $P\cap Z=Q\cap Z$. Then{\rm:}

\item(i)
$I_r(M/MK)$ coincides with the image of $I_r(M)$ in $R/K$.

\item(ii)
$I_r(M)\sbs K$ if and only if $(M/MK)^l\cong(R/K)^n$ in $\M_R$.

\item(iii)
If $I_r(M)\sbs K$ then $r_P(M)=r$ for each $P\in\Max R$ with $P\sps K$.

\endproclaim

\Proof.
If $\m=P\cap Z$ where $P\in\Max R$, $\,P\sps K$, then we have $r_Q(M)\le r$ 
for all $Q\in\Max_\m R$ by the hypothesis; hence $\m\in U_r(M)$ according to 
Lemma 2.3. The Remark following Lemma 1.10 now proves (i). Since $r_P(M)l\le 
n$ for each $P\in\Max R$ with $P\sps K$, the $R/K\!$-module $(M/MK)^l$ is 
$n$-generated by Lemma 1.6. Hence (ii) follows from (i) and Lemma 1.4(ii).
If $(M/MK)^l\cong(R/K)^n$, then $(M/MP)^l\cong(R/P)^n$ for any $P\in\Max R$ 
with $P\sps K$; the comparison of lengths of the two modules appearing in the 
latter isomorphism yields (iii).
\endproof

\proclaim
Lemma 2.6.
Suppose that there is an integer $n\ge0$ such that $M_\m\cong R^n_\m$ in 
$\M_{R_\m}$ for each $\m\in\Max Z$. Then $M$ is projective{\rm;} $M$ is a 
generator in $\M_R$ provided $M\ne0$. For each $\m\in\Max Z$ there exists 
$z\in Z\setm\m$ such that $M_z\cong R^n_z$ in $\M_{R_z}$.
\endproclaim

\Proof.
We first prove that $M$ is finitely presented. Consider any $\M_R$-epimorphism
$\ph:R^m\to M$ with kernel $K$. We have to check that $K$ is finitely 
generated. For each $\m\in\Max Z$ the localization $\ph_\m$ of $\ph$ at $\m$ 
gives rise to an exact sequence of $R_\m$-modules $0\to K_\m\to R^m_\m\to 
M_\m\to0$. Since $M_\m$ is free, the sequence splits, and so $R^m_\m\cong 
K_\m\oplus R^n_\m$. It follows that $m\ge n$ and $K_\m\cong R^{m-n}_\m$ by the 
cancellation property for projective modules over a semilocal ring \cite{Bass, 
Ch.\ IV, (1.4)}.

For any fixed $\m$ we can pick $m$ elements $v_1,\ldots,v_m\in R^m$ with the 
property that $v_{n+1},\ldots,v_m\in K$, the images of $v_{n+1},\ldots,v_m$ in 
$K_\m$ give a basis for $K_\m$ over $R_\m$, and the images of $v_1,\ldots,v_n$ 
in $R^m_\m$ give a basis for a complementary summand. Let $K'$ and $N$ be the 
submodules of $R^m$ generated by $v_{n+1},\ldots,v_m$ and $v_1,\ldots,v_n$, 
respectively. Since $R^m_\m=K'_\m+N_\m$, there exists $z\in Z\setm\m$ such 
that $R^mz\sbs K'+N$. For any $\n\in\Max Z$ with $z\notin\n$ we have 
$R^m_\n=K'_\n+N_\n$. Since $K'\sbs K$, it follows that $\ph_\n(N_\n)=M_\n$. 
The $R_\n$-module $N_\n$ is generated by $n$ elements. Since the ring $R_\n$ 
is weakly finite and $M_\n\cong R^n_\n$, the images of those elements in 
$M_\n$ are a basis for $M_\n$ over $R_\n$. In other words, $\ph_\n$ induces an 
isomorphism of $N_\n$ onto $M_\n$. Hence $K_\n\cap N_\n=0$, and therefore 
$K_\n=K'_\n$.

Denote by $\U$ the collection of all open subsets $U$ of $\Max Z$ with the 
property that there exists a finitely generated submodule $L\sbs K$, depending 
on $U$, such that $K_\n=L_\n$ for all $\n\in U$. We have just proved that each
$\m\in\Max Z$ has an open neighborhood contained in $\U$. It is also clear 
that $U\cup U'\in\U$ whenever $U,U'\in\U$. Since the space $\Max Z$ is 
quasicompact, we conclude that $\Max Z\in\U$. This means that there exists a 
finitely generated submodule $L\sbs K$ such that $K_\n=L_\n$ for all 
$\n\in\Max Z$. But then $K=L$. Thus $M$ is finitely presented, as claimed.

For any fixed $\m$ we can find an $\M_R$-morphism $\psi:R^n\to M$ whose 
localization $R^n_\m\to M_\m$ is an isomorphism. Since $M$ is finitely 
generated, there exists $s\in Z\setm\m$ such that the $\M_{R_s}$-morphism
$\psi_s:R^n_s\to M_s$ induced by $\psi$ is surjective. Since the $R_s$-module 
$M_s$ is finitely presented, $\Ker\psi_s$ is finitely generated. Then
$\Ker\psi_s$ is annihilated by some element $z\in Z\setm\m$. We may assume 
that $Z_z$ is a localization of $Z_s$, in which case $\psi_s$ induces an 
isomorphism $R^n_z\to M_z$.

For $V\in\M_R$ and $\m\in\Max Z$ the canonical map
$$
\Hom_R(M,V)\ot_ZZ_\m\to\Hom_{R_\m}(M_\m,V_\m)
$$
is bijective by \cite{Bass, Ch.\ III, (4.5)}. Let now 
$\xi:\Hom_R(M,V)\to\Hom_R(M,W)$ be the map induced by an $\M_R$-epimorphism 
$V\to W$. Since $M_\m$ is projective in $\M_{R_\m}$, the map
$\xi\ot\id:\Hom_R(M,V)\ot_ZZ_\m\to\Hom_R(M,W)\ot_ZZ_\m$ is surjective for 
each $\m$, but then $\xi$ is itself surjective. This proves that $M$ is
projective in $\M_R$.

Suppose that $M\ne0$. Then $n\ne0$, and so $M_\m\ne0$ for each
$\m\in\Max Z$. Denote by $T$ the trace ideal of $M$. Thus $T=\sum\,f(M)$
where $f$ runs over $\Hom_R(M,R)$. Since 
$\Hom_{R_\m}(M_\m,R_\m)\cong\Hom_R(M,R)\ot_ZZ_\m$,
the trace ideal of the $R_\m$-module $M_\m$ coincides with $T_\m$. It follows 
that $T_\m=R_\m$ since $M_\m$ is free. As this is valid for each $\m$,
we get $T=R$. This means that $M$ is a generator.
\endproof

\proclaim
Lemma 2.7.
Let $P\in\Spec R$ and $\m\in\Max Z$. If $P\cap Z\sbs\m$ then $P\sbs Q$ for 
some $Q\in\Max_\m R$.
\endproclaim

\Proof.
The hypothesis implies that $1\notin P_\m$. So $P_\m$ is contained in a 
maximal ideal of $R_\m$, that is, an ideal $Q_\m$ for some $Q\in\Max_\m R$. 
Then $P\sbs Q$.
\endproof

Recall that a {\it Jacobson ring} is a ring in which every prime ideal is an 
intersection of primitive ideals.

\proclaim
Lemma 2.8.
Let $P\in\Spec R$ and $\p=P\cap Z$. Suppose that $Z$ is a Jacobson ring and
$R_\p/\p R_\p$ is an artinian ring. Then for any $z\in Z\setm\p$ the 
intersection $K$ of all ideals $Q\in\Max R$ such that $P\sbs Q$ and $z\notin Q$
coincides with $P$. 
\endproclaim

\Proof.
Clearly $P\sbs K$. In view of Lemma 2.7 $K\cap Z$ coincides with the 
intersection $\b$ of all ideals $\m\in\Max Z$ such that $\p\sbs\m$ and
$z\notin\m$. Since $Z$ is a Jacobson ring, we have $\p=\a\cap\b$ where $\a$ 
denotes the intersection of all ideals $\m\in\Max Z$ such that $\p\sbs\m$ and 
$z\in\m$. In particular, $\a\b\sbs\p$. Since $z\in\a$, we have $\a\not\sbs\p$; 
hence $\b\sbs\p$ because $\p$ is prime. We conclude that $K\cap Z=\p$.

Let $R'=R_\p/P_\p\cong R/P\ot_ZZ_\p$. Since $R/P$ is a prime ring, its central 
subring $Z/\p$ contains no zero divisors of $R/P$ other than $0$. Hence 
$R/P$ is embedded into $R'$. The ring $R'$, as a homomorphic image of 
$R_\p/\p R_\p$, is artinian. Each nonzero ideal of $R'$ intersects $R/P$ 
nontrivially. It follows that $R'$ is prime, in which case $R'$ is actually 
simple. If the ideal $K'$ of $R'$ generated by the image of $K$ contained $1$,
$K/P$ would have a nonzero intersection with $Z/\p$, which is impossible.
We must have $K'=0$, which entails $K/P=0$, i.e., $K=P$.
\endproof

\Remark.
If $R$ is module-finite over $Z$, then $R_\p/\p R_\p$ is a finite dimensional 
algebra over a field, so that the artinian hypothesis in Lemma 2.8 is 
fulfilled. In this case $\p R_\p\sbs\Jac(R_\p)$ by \cite{Az51, Corollary to 
Lemma 2} or \cite{Cur53, Lemma 3.1}, which implies that $R_\p$ is semilocal.

It is well-known that the Jacobson property goes up from $Z$ to $R$ in 
the module-finite case. The first result of this kind, due to Curtis 
\cite{Cur53, Th. 4.3}, assumed ACC on $Z$-submodules of $R$. Subsequently 
several generalizations have been found, e.g. \cite{Cor88}, \cite{Pro67}, 
\cite{Rob81}. Under previous assumptions $R_z$ is module-finite over $Z_z$ and 
$Z_z$ is Jacobson by \cite{Bou, Ch.~V, \S3, Th.~3}; so $R_z$ is Jacobson for 
any $z\in Z$. This reduces to the conclusion of Lemma 2.8.
\endremark

\section
3. Semilocal factor algebras of module algebras

Suppose that $H$ is a Hopf algebra and $A$ a left $H$-module algebra over the 
ground field $k$. The compatibility of the $H$-module structure with the 
algebra structure on $A$ is expressed by means of the identities
$$
h1_A=\ep(h)1_A,\qquad h(ab)=\sum\,\,(h\1a)(h\2b)
$$
where $h\in H$, $a,b\in A$, and $1_A$ is the unity element of $A$. For an 
ideal $I$ of $A$ and a subcoalgebra $C$ of $H$ put
$$
I_C=\{a\in A\mid Ca\sbs I\}.
$$
Clearly $I_C$ is also an ideal of $A$. In particular, $I_H$ is the largest 
$H$-stable ideal of $A$ contained in $I$. If $C,C'$ are two subcoalgebras with 
$C\sbs C'$ then $I_C\sps I_{C'}$.

Our subsequent arguments require the factor algebras $A/I_C$ to be semilocal. 
We wish to know the cases in which this property of $A/I_C$ can be 
established.

We always consider $\Hom(C,A/I)$ equipped with the convolution multiplication.
For $a\in A$ define $\at\in\Hom(C,A/I)$ by the rule
$$
\at(c)=ca+I,\qquad c\in C.
$$
The map $\tau:A\to\Hom(C,A/I)$ given by the assignment $a\mapsto\at$ is a
homomorphism of algebras and $\Ker\tau=I_C$. The inclusion
$k\hookrightarrow A/I$ allows us to identify the dual algebra $C^*$ of $C$
with a subalgebra of $\Hom(C,A/I)$.

\proclaim
Lemma 3.1.
If $\dim C<\infty$, then $\Hom(C,A/I)=\tau(A)C^*$ and so $\Hom(C,A/I)$ is 
left module-finite over $\tau(A)$.
\endproclaim

\Proof.
The verification is straightforward. In case $I=0$ the statement is contained 
in \cite{Sk-Oy07, Lemma 2.1(iv)}. The general case is immediate since the 
canonical projection $A\to A/I$ induces a surjective algebra homomorphism 
$\Hom(C,A)\to\Hom(C,A/I)$.
\endproof

\Remark.
It is probably not true that $\Hom(C,A/I)$ is right module-finite over 
$\tau(A)$ in case of Hopf algebras whose antipode is not bijective. This is 
essentially the reason for our use of left side conditions in this section.
\endremark

\proclaim
Lemma 3.2.
If $C$ and $A/I$ are finite dimensional, then so too is $A/I_C$.
\endproclaim

This is clear since $A/I_C$ is embedded into $\Hom(C,A/I)$.

\proclaim
Lemma 3.3.
Let $T$ be a ring left module-finite over a subring $R$. Suppose that $R$ is 
either {\rm(a)} left noetherian or {\rm(b)} left module-finite over a 
commutative subring $R'$. If $T$ is left artinian {\rm(}semiprimary, 
semilocal{\rm),} then so too is $R$.
\endproclaim

\Proof.
(i) If $T$ is left artinian then $R$ is left artinian by Bj\"ork's results.
In case (a) \cite{Bj73, Cor.\ 0.2} applies. In case (b) $T$ is left 
module-finite over $R'$; so $R'$ is artinian by \cite{Bj73, Th.\ 3.3}. Then 
$T$ has finite length as an $R'$-module with respect to left multiplications, 
and the same holds for $R$.

(ii) Suppose that $T$ is semiprimary and $J=\Jac(T)$. Part (i) shows that the 
subring $R/(J\cap R)$ of the artinian ring $T/J$ is left artinian. Since 
$J\cap R$ is a nilpotent ideal of $R$, it is clear that $R$ is semiprimary.

(iii) A result of Camps and Dicks \cite{Camps93} says that a subring of a 
semilocal ring is itself semilocal provided that the subring is {\it full}, 
that is, each non-invertible element of the subring is not invertible in the 
ambient ring. We will check that $R$ is a full subring of $T$; it will follow 
then that $R$ is semilocal whenever so is $T$. Let $x\in R$ be invertible in 
$T$. We have to show that $x^{-1}\in R$. In case (a) $T$ is a noetherian 
$R$-module on the left side. Hence the chain of submodules 
$R\sbs Rx^{-1}\sbs Rx^{-2}\sbs\cdots$ is ultimately constant, i.e.,
$x^{-n}\in Rx^{1-n}$ for some $n>0$. Multiplying by $x^{n-1}$ proves the claim.

In case (b) $T$ is a finitely generated $R'$-module on the left side. The 
right multiplication by $x^{-1}$ defines an endomorphism $f$ of that module. 
Since $R'$ is commutative, $f$ satisfies an equation 
$f^n=\sum_{i=0}^{n-1}c_if^i$ for some $n>0$ and $c_0,\ldots,c_{n-1}\in R'$. 
Then we have $tx^{-n}=\sum_{i=0}^{n-1}c_itx^{-i}$ for all $t\in T$. 
Substituting $t=1$, we deduce that $x^{-n}\in\sum_{i=0}^{n-1}R'x^{-i}\sbs 
Rx^{1-n}$, which leads to the desired conclusion.
\endproof

\proclaim
Lemma 3.4.
Suppose that $\dim C<\infty$ and $A$ is either left noetherian or left 
module-finite over a commutative subring. If $A/I$ is left artinian 
{\rm(}semiprimary, semilocal{\rm),} then so too is $A/I_C$.
\endproclaim

\Proof.
We take $T=\Hom(C,A/I)$ and $R=\tau(A)$. By Lemma 3.1 $T$ is left module-finite
over $R$. Let us identify $T$ with the algebra $C^*\ot A/I$ by means of 
the canonical isomorphism. Thus $T$ is left module-finite over the subring 
$1\ot A/I$ isomorphic to $A/I$. If $A/I$ is left artinian, so is $T$.
Since the finite dimensional subalgebra $C^*\ot1$ centralizers $1\ot A/I$, 
the ideal $J=C^*\ot\Jac(A/I)$ of $T$ is contained in the Jacobson radical of 
$T$ \cite{Lam, Prop. 5.7}. If $A/I$ is semilocal, then $T/J$ is artinian, and 
it follows that $T$ is semilocal. If $A/I$ is semiprimary then $J$ is 
nilpotent; hence $T$ is semiprimary. Since $A/I_C\cong\tau(A)$, an application 
of Lemma 3.3 yields all conclusions of Lemma 3.4.
\endproof

\Remark.
The conclusion of Lemma 3.4 can be rephrased in the language of \cite{Mo95, 
Def. 3.4} as follows: the action of $H$ on $A$ is $\F$-continuous where $\F$ 
is the filter consisting of those ideals $I$ of $A$ for which $A/I$ is left 
artinian in one case, semiprimary in the second and semilocal in the third. 
Another result of this kind will be presented in Lemma 3.6.
\endremark

\proclaim
Lemma 3.5.
Let $Z\sbs R\sbs T$ be a tower of rings where $Z$ is central in $R$, the ring 
$R$ has semilocal localizations with respect to $Z$ and $T$ is left
module-finite over $R$. If $T$ is semilocal then $R$ is semilocal. If $T$ is 
semiprimary then so too is $Z${\rm;} if also all rings $R_\m/\m R_\m$ with 
$\m\in\Max Z$ are semiprimary then $R$ is semiprimary.
\endproclaim

\Proof.
Note that the version of Lemma 2.2 for left $R$-modules is also valid since 
we may replace $R$ with the opposite ring. Take $M=T$ regarded as an 
$R$-module with respect to left multiplications. Lemma 2.2 shows that the 
equality $aT=T$ for $a\in Z$ implies $aZ=Z$. In other words, $a^{-1}\in Z$ 
whenever $a$ is invertible in $T$. So $Z$ is a full subring of $T$, and the 
Camps-Dicks Theorem ensures that $Z$ is semilocal (cf.\ the proof of Lemma 
3.3). Since the set $\Max Z$ is finite, by Lemma 2.1 $R$ has finitely many 
right primitive factor algebras $R/P$, and each of those is simple artinian. 
Hence $R$ is semilocal.

Suppose that $T$ is semiprimary. Since $J=\Jac(T)$ is nilpotent, the ring $R$
(resp., $Z$) is semiprimary if and only if so is $R/(R\cap J)$ (resp.,
$Z/(Z\cap J)$). Passing to the tower of rings
$Z/(Z\cap J)\sbs R/(R\cap J)\sbs T/J$, we may assume that $T$ is artinian. 
Given $a\in Z$, there exists an integer $n>0$ such that $a^nT=a^{n+1}T$. Now 
$a^nT$ is a finitely generated left $R$-submodule of $T$ since $aR=Ra$. 
Applying Lemma 2.2 with $M=a^nT$ and $\a=aZ$, we deduce that $(1-b)a^nT=0$ for 
some $b\in aZ$. If $a\in\Jac(Z)$, then $1-b$ is invertible, whence $a^nT=0$, 
i.e., $a^n=0$. This shows that $\Jac(Z)$ is nil. Since nil subrings of 
artinian rings are nilpotent, $\Jac(Z)$ is nilpotent. This means that $Z$ is 
semiprimary. Since $Z$ is commutative, $Z$ is the finite direct product of 
local rings $Z_\m$, $\,\m\in\Max Z$, with nilpotent maximal ideals $\m Z_\m$. 
Then $R\cong\prod R_\m$ and $\m R_\m$ is a nilpotent ideal of $R_\m$ for each 
$\m$. If $R_\m/\m R_\m$ is semiprimary, so too is $R_\m$. When all rings 
$R_\m$ are semiprimary, $R$ is semiprimary.
\endproof

\setitemsize(iii)
\proclaim
Lemma 3.6.
Suppose that $\dim C<\infty$ and $A$ has semilocal localizations with 
respect to a central subring $Z$.

\item(i)
If $A/I$ is semilocal then $A/I_C$ is semilocal.

\item(ii)
If $A/I$ is semiprimary then $Z/(Z\cap I_C)$ is semiprimary.

\item(iii)
If $A/I$ and all rings $A_\m/\m A_\m$ are semiprimary then $A/I_C$ is 
semiprimary.

\endproclaim

\Proof.
We apply Lemma 3.5 to the tower $\tau(Z)\sbs\tau(A)\sbs\Hom(C,A/I)$. As
pointed out in the proof of Lemma 3.4, $\Hom(C,A/I)$ is semilocal 
(semiprimary) whenever so is $A/I$.
\endproof

\section
4. The orbit relation on the maximal spectrum

We continue to assume that $A$ is an $H$-module algebra. For $P,Q\in\Max A$ 
define $P\ds_HQ$ if $P_C\sbs Q$ for some finite dimensional subcoalgebra 
$C\sbs H$. Here $P_C$ denotes the ideal of $A$ defined in section 3.

\proclaim
Lemma 4.1.
The relation $\ds_H$ is reflexive and transitive.
\endproclaim

\Proof.
If $C=k$, then $P_C=P$. Since $P\sbs P$, we get $P\ds_HP$. Suppose that 
$P,P',P''\in\Max A$ satisfy $P\ds_HP'$ and $P'\ds_HP''$. Then $P_C\sbs P'$ and 
$P'_{C'}\sbs P''$ for some finite dimensional subcoalgebras $C,C'\sbs H$. Note 
that $CC'$ is also a finite dimensional subcoalgebra of $H$. If $a\in 
P_{CC'}$, that is, $CC'a\sbs P$, then $C'a\sbs P_C\sbs P'$, whence $a\in 
P'_{C'}\sbs P''$. This shows that $P_{CC'}\sbs P''$, and therefore 
$P\ds_HP''$.
\endproof

If $H=kG$ is a group algebra, then any finite dimensional subcoalgebra $C$ of 
$H$ is spanned by a finite subset, say $X$, of $G$. Clearly
$P_C=\bigcap_{g\in X}g^{-1}(P)$. If $Q\in\Max A$ contains $P_C$, then $Q$ 
contains the product of the ideals $g^{-1}(P)$, $\,g\in X$, taken in any 
order; since $Q$ is prime, $Q\sps g^{-1}(P)$ for some $g\in X$. The maximality 
of $P$ ensures then that $Q=g^{-1}(P)$. Thus $P\ds_HQ$ if and only if $P$ and 
$Q$ lie in the same $G$-orbit.

The previous example suggests that $\ds_H$ may also be symmetric, that is, an 
equivalence relation on $\Max A$ in general. It is not clear whether this is 
always true. We will be able to provide a confirmation in several cases. When 
the relation $\ds_H$ is symmetric, we call it the {\it$H\!$-orbit 
equivalence relation}.

Note that $P_H$ coincides with the intersection of the family of ideals $P_C$
with $C$ a finite dimensional subcoalgebra. It follows that $P_H\sbs Q_H$ 
whenever $P\ds_HQ$. If $P\ds_HQ$ and $Q\ds_HP$ then $P_H=Q_H$, that is, $P$ 
and $Q$ belong to the same $H\!$-stratum, in the language of \cite{Br}. In 
general the $H\!$-stratification defines a coarser equivalence relation.

The proof of the next lemma uses essentially the same argument as given by 
Chin \cite{Ch90, Lemma 2.2} in the case where $H$ is finite dimensional and 
pointed; it was further generalized by Montgomery and Schneider \cite{Mo99, 
Th. 3.7}.

\proclaim
Lemma 4.2.
Suppose that $H'$ is a Hopf subalgebra of $H$ containing the coradical of $H$.
Then for $P,Q\in\Max A$ one has $P\ds_HQ$ if and only if $P\ds_{H'}Q$.
The relation $\ds_H$ is symmetric if and only if so is $\ds_{H'}$.
\endproclaim

\Proof.
Suppose that $C$ is a finite dimensional subcoalgebra of $H$ and $C_0$ 
denotes the coradical of $C$. Consider the coradical filtration $C_0\sbs 
C_1\sbs\cdots$ of $C$. As $\dim C<\infty$, we have $C_n=C$ for some $n$. Let 
$P\in\Max A$. We will prove by induction on $i\ge0$ that $P_{C_0}^{i+1}\sbs
P_{C_i}$. For $i=0$ this is clear. Suppose that the claim is valid for some 
$i\ge0$ and $c\in C_{i+1}$. Since $\De(c)\in C_0\ot C+C\ot C_i$, we deduce
$$
c(P_{C_0}^{i+2})\sbs\sum_{(c)}c_{(1)}(P_{C_0})\cdot 
c_{(2)}(P_{C_0}^{i+1})\sbs P,
$$
showing that $P_{C_0}^{i+2}\sbs P_{C_{i+1}}$. In particular, 
$P_{C_0}^{n+1}\sbs P_C$. It follows that for $Q\in\Max A$ the inclusions 
$P_C\sbs Q$ and $P_{C_0}\sbs Q$ are equivalent to each other. Since $C_0\sbs 
H'$, we conclude that $P\ds_HQ$ if and only if $P\ds_{H'}Q$.
\endproof

\proclaim
Corollary 4.3.
If $H$ is pointed with the group $G$ of grouplike elements then $P\ds_HQ$ for 
$P,Q\in\Max A$ if and only if $P$ and $Q$ lie in the same $G$-orbit. 
\endproclaim

\Proof.
In this case the coradical of $H$ coincides with the group algebra $kG$.
\endproof

\proclaim
Proposition 4.4.
If $A$ is right artinian then $\ds_H$ is symmetric.
\endproclaim

\Proof.
If $C$ and $C'$ are two finite dimensional subcoalgebras of $H$, so also 
is $C+C'$, and $P_{C+C'}=P_C\cap P_{C'}$. Since $A$ satisfies DCC on right 
ideals, the set of ideals $P_C$ with $C$ a finite dimensional subcoalgebra of 
$H$, contains a smallest element which has to coincide with $P_H$. Hence for 
$P,Q\in\Max A$ one has $P\ds_HQ$ if and only if $P_H\sbs Q$. The right artinian
$H$-module algebra $A/P_H$ has a maximal ideal $P/P_H$ which contains no 
nonzero $H$-stable ideals of $A/P_H$. By \cite{Sk-Oy07, Lemma 4.2} $A/P_H$ is 
$H$-simple. If $P_H\sbs Q$, then $Q_H$ is an $H$-stable ideal of $A$ 
containing $P_H$, and we must have $Q_H=P_H$. The inclusion $Q_H\sbs P$ 
entails $Q\ds_HP$.
\endproof

\proclaim
Proposition 4.5.
Suppose that $X\sbs\Max A$ is a subset such that for each $P\in X$ and 
each finite dimensional subcoalgebra $C$ of $H$ the factor ring $A/P_C$ is 
semiprimary and each maximal ideal of $A$ containing $P_C$ lies in $X$. If 
either {\rm(a)} $\dim H<\infty$ or {\rm(b)} $H$ is generated by a family $\H$ 
of Hopf subalgebras such that the relation $\ds_{H'}$ is symmetric on $X$ for 
each $H'\in\H$, then the relation $\ds_H$ is symmetric on $X$.
\endproclaim

\Proof.
In case (a) $P_H$ is the smallest element in the set of ideals $P_C$ with $C$ 
a finite dimensional subcoalgebra. Hence for $P,Q\in\Max A$ one has $P\ds_HQ$ 
if and only if $P_H\sbs Q$. The $H$-module algebra $A/P_H$ is semiprimary by 
the hypothesis. Its maximal ideal $P/P_H$ contains no nonzero $H$-stable 
ideals of $A/P_H$. It follows that $A/P_H$ is $H$-semiprime, i.e., $A/P_H$ has 
no nonzero $H$-stable nilpotent ideals. By \cite{Sk-Oy07, Th.\ 0.3 and Lemma 
4.2} $A/P_H$ is $H$-simple. Then the inclusion $P_H\sbs Q$ implies $P_H=Q_H$, 
and so $Q\ds_HP$.

Assume now that $H$ satisfies condition (b). Denote by $\C$ the collection of 
subcoalgebras $C$ of $H$ such that $\dim C<\infty$ and for any pair $P,Q\in X$ 
satisfying $P_C\sbs Q$ one has $Q\ds_HP$. By the hypothesis $\C$ contains all
finite dimensional subcoalgebras of any $H'\in\H$.

We claim that $C+C'\in\C$ and $CC'\in\C$ whenever $C$, $C'$ are both from 
$\C$. Suppose that $P,Q\in X$ are such that $P_{C+C'}\sbs Q$. Since
$P_{C+C'}=P_C\cap P_{C'}\sps P_CP_{C'}$ and $Q$ is a prime ideal, we have 
either $P_C\sbs Q$ or $P_{C'}\sbs Q$, whence $Q\ds_HP$. This proves the first
inclusion in our claim.

We also have to show that $Q\ds_HP$ whenever $P,Q\in X$ satisfy $P_{CC'}\sbs 
Q$. Denote by $Y$ the set of maximal ideals of $A$ containing $P_C$. By the 
hypothesis $Y\sbs X$, and $Y$ is finite since $A/P_C$ is semiprimary. If $J$ 
denotes the intersection of all ideals from $Y$, then $J/P_C$ coincides with 
the Jacobson radical of $A/P_C$, which is nilpotent. It follows that there 
exists a finite sequence $Q_1,\ldots,Q_n$ of ideals from $Y$ (with repetitions 
allowed) such that $Q_1\cdots Q_n\sbs P_C$. If $a_1,\ldots,a_n\in A$ are any 
elements such that $C'a_i\sbs Q_i$ for each $i=1,\ldots,n$, then $C'(a_1\cdots 
a_n)\sbs Q_1\cdots Q_n\sbs P_C$, and then $CC'(a_1\cdots a_n)\sbs P$. This 
shows that $(Q_1)_{C'}\cdots(Q_n)_{C'}\sbs P_{CC'}\sbs Q$. Since $Q$ is prime,
we must have $(Q_i)_{C'}\sbs Q$ for at least one $i$. The inclusions 
$C',C\in\C$ imply that $Q\ds_HQ_i$ and $Q_i\ds_HP$. The transitivity of the 
relation $\ds_H$ entails $Q\ds_HP$, as required.

It is clear now that the union $U$ of all coalgebras from $\C$ is a subalgebra 
of $H$. If $H'\in\H$ then $H'\sbs U$ since $H'$ is the union of its finite 
dimensional subcoalgebras. Since $H$ is generated by $\H$, we get $U=H$.
Each finite dimensional subcoalgebra $C$ of $H$ is contained therefore in some 
$C'\in\C$; since $P_C\sps P_{C'}$ for any $P\in X$, it is clear that $C\in\C$.
Thus $Q\ds_HP$ whenever $P,Q\in X$ satisfy $P_C\sbs Q$.
\endproof

There is a different interpretation of the relation $P\ds_HQ$ in terms
of certain operations with modules. Denote by $\M^H$ the category of right 
$H$-comodules. Given $U\in\M^H$ and $V\in\M_A$, we define right $A$-module 
structures on vector spaces $U\ot V$ and $\Hom(U,V)$ by the rules
$$
\eqalign{
(u\ot v)a&{}=\sum\,u\0\ot v\bigl((Su\1)a\bigr),\cr
(\eta a)(u)&{}=\sum\,\eta(u\0)(u\1a)\cr
}
$$
where $u\in U$, $v\in V$, $a\in A$, $\eta\in\Hom(U,V)$ and $S:H\to H$ is the 
antipode (see \cite{Sk-Oy07, section 1}). If $\dim U<\infty$ then $U^*$ is a 
right $H$-comodule with structure map $U^*\to U^*\ot H$, 
$\,\xi\mapsto\sum\xi\0\ot\xi\1$, such that
$$
\sum\,\xi\0(u)\xi\1=\sum\,\xi(u\0)Su\1
$$
for all $u\in U$. Note that the evaluation map $\ev:U^*\ot U\to k$ is an 
$\M^H$-morphism provided $k$ has the trivial comodule structure.

\setitemsize(ii)
\proclaim
Lemma 4.6.
Let $U\in\M^H$ and $V, W\in\M_A$.

\item(i)
$\Hom_A(U\ot V,W)\cong\Hom_A\bigl(V,\,\Hom(U,W)\bigr)$.

\item(ii)
If\/ $\dim U<\infty$ then $U\ot V\cong\Hom(U^*\!,V)$ in $\M_A$.

\endproclaim

\Proof.
(i) This is the isomorphism from \cite{Sk-Oy07, Lemma 1.1}. It is induced by the
canonical linear bijection $\Hom(U\ot V,W)\cong\Hom\bigl(V,\,\Hom(U,W)\bigr)$.

(ii) We obtain
$\Hom_A\bigl(U^*\ot(U\ot V),V\bigr)\cong\Hom_A\bigl(U\ot V,\,\Hom(U^*\!,V)\bigr)$
as a special case of (i). The canonical map $\ph:U\ot V\to\Hom(U^*\!,V)$ 
corresponds to the composite
$$
U^*\ot(U\ot V)\cong(U^*\ot U)\ot V\lmapr4{\ev\ot\id}k\ot V\cong V.
$$
Since the latter is an $\M_A$-morphism by functoriality, so too is $\ph$. The 
assumption $\dim U<\infty$ entails the bijectivity of $\ph$. Thus $\ph$ is an 
isomorphism in $\M_A$.
\endproof

\setitemsize(iii)
\proclaim
Lemma 4.7.
Suppose that $A/P$ and $A/Q$ are simple artinian. Let $V$ and $W$ be simple 
right $A$-modules whose annihilators coincide with $P$ and $Q$, respectively. 

\item(i)
$P\ds_HQ$ if and only if $W$ is a subfactor of the right $A$-module 
$\Hom(U,V)$ for some finite dimensional $U\in\M^H$.

\item(ii)
If $W$ is a submodule of $\Hom(U,V)$ then both $P\ds_HQ$ and $Q\ds_HP$ hold.

\item(iii)
If $W$ is a factor module of $\Hom(U,V)$ and the antipode of $H$ is bijective 
then $P\ds_HQ$ and $Q\ds_HP$ too.

\endproclaim

\Proof.
(i) There is an isomorphism $A/P\cong V^n$ in $\M_A$ for some integer $n>0$. 
Suppose that there exists a finite dimensional subcoalgebra $C$ of $H$ such 
that $P_C\sbs Q$. We may regard $C$ as a right $H$-comodule with respect to 
the comultiplication. The right $A$-module structure on $\Hom(C,A/P)$ derives 
from the algebra homomorphism $\tau:A\to\Hom(C,A/P)$ defined in section 3. 
Since $\Ker\tau=P_C$, the factor algebra $A/P_C$ is embedded in $\Hom(C,A/P)$.  
As $W$ is a simple $A/P_C$-module, $W$ is a subfactor of $\Hom(C,A/P)$ as a 
right $A$-module. The latter module is the direct sum of $n$ copies of 
$\Hom(C,V)$. Hence $W$ is a subfactor of $\Hom(C,V)$.

Conversely, suppose that $W$ is a subfactor of $\Hom(U,V)$ for some finite 
dimensional $U\in\M^H$. Since $H$ is an injective cogenerator in $\M^H$, there 
exists a monomorphism $\ph:U\to H^m$ in $\M^H$ for some integer $m>0$. Then 
$\ph(U)\sbs C^m$ for a suitable finite dimensional subcoalgebra $C\sbs H$. 
Hence $\Hom(U,V)$ is a homomorphic image of the right $A$-module 
$\Hom(C^m,V)$. It follows that $W$ is a subfactor of $\Hom(C,V)$, and also of 
$\Hom(C,A/P)$. Using again the equality $\Ker\tau=P_C$ from the previous 
paragraph, we deduce that $P_C$ annihilates $W$, whence $P_C\sbs Q$.

(ii) If there exists an $\M_A$-monomorphism $W\to\Hom(U,V)$ then there also 
exists a nonzero $\M_A$-morphism $U\ot W\to V$ by Lemma 4.6(i). The latter 
has to be surjective since $V$ is simple. Lemma 4.6(ii) shows that $V$ is a 
subfactor of $\Hom(U^*\!,W)$, whence $Q\ds_HP$ by part (i).

(iii) Suppose that the antipode is bijective. Then any $U\in\M^H$, $\,\dim 
U<\infty$, is isomorphic to $(U')^*$ for some finite dimensional 
$U'\in\M^H$. By Lemma 4.6(ii) $\Hom(U,V)\cong U'\ot V$. If there exists an 
$\M_A$-epimorphism $\Hom(U,V)\to W$ then there also exists a nonzero 
$\M_A$-morphism $V\to\Hom(U',W)$ by Lemma 4.6(i). In this case $V$ is 
isomorphic with a submodule of $\Hom(U',W)$, whence $Q\ds_HP$.
\endproof

\Remark.
We may regard $\M_A$ as a right module category over the tensor category
$(\M^H)\op$, opposite to $\M^H$, with respect to the bifunctor 
$(V,U)\mapsto\Hom(U,V)$. Lemma 4.7 shows that $\ds_H$ corresponds to a certain
relation on the set of isomorphism classes of simple right $A$-modules defined 
in purely categorical terms. In case of an arbitrary left module category $\M$
over a finite tensor category $\C$ such a relation was introduced by Etingof 
and Ostrik \cite{Et04, Lemma 3.8}. It was proved there that this relation is 
symmetric under the assumption that $\C$ has projective covers and $P\ot X$ is
projective in $\M$ for any projective object $P\in\C$ and any object $X\in\M$.
The second condition is rather nontrivial to verify.
\endremark

\proclaim
Lemma 4.8.
Suppose that $A$ has semilocal localizations with respect to a central 
subring $Z$. Given $P,Q\in\Max A$ with $P\ds_HQ$, let $\n=Q\cap Z$. If either 
{\rm(a)} $A_\n/\n A_\n$ is semiprimary or {\rm(b)} the antipode of $H$ is 
bijective{\rm,} then there exists $Q'\in\Max_\n A$ satisfying $P\ds_HQ'$ and 
$Q'\ds_HP$.
\endproclaim

\Proof.
Recall from Lemma 2.1 that both rings $A/P$ and $A/Q$ are simple artinian.  
Let $V$ and $W$ be as in Lemma 4.7. Then $W$ is a subfactor of the right 
$A$-module $M=\Hom(U,V)$ for some finite dimensional $U\in\M^H$. Denoting by 
$\rho:U\to U\ot H$ the comodule structure map, we have $\rho(U)\sbs U\ot C$ 
for some finite dimensional subcoalgebra $C\sbs H$. Since $VP=0$, it is 
immediate from the definition of the $A$-module structure on $M$ that $P_C$ 
annihilates $M$. Put $\a=P_C\cap Z$. The commutative ring $Z/\a$ is 
semiprimary by Lemma 3.6; it is therefore a finite direct product of primary 
rings. Then $A/\a A\cong\prod_{\m\in X}A_\m/\a A_\m$ where $X$ is the finite 
set of those $\m\in\Max Z$ for which $\a\sbs\m$. Since $M\a=0$, we have 
$M\cong\prod_{\m\in X}M_\m$, and $A$ operates in $M_\m$ via the projection 
onto $A_\m/\a A_\m$. Then $W$ is a subfactor of $M_\m$ for some $\m\in X$. 
Since $W\n=0$, while all elements of $Z\setm\m$ are invertible on $W$, we must 
have $\n\sbs\m$. As $\n\in\Max Z$ by Lemma 2.1, this yields $\n=\m$, showing 
that $\n\in X$ and $M_\n\ne0$.

Suppose that (a) holds. The primary ring $Z_\n/\a Z_\n$ has a nilpotent 
maximal ideal generated by $\n$. Hence $\n A_\n/\a A_\n$ is a nilpotent ideal 
of $A_\n/\a A_\n$, and it follows that $A_\n/\a A_\n$ is semiprimary. Since 
$M_\n$ is a nonzero right $A_\n/\a A_\n$-module, it contains a simple 
submodule, say $W'$. Denote by $Q'$ the annihilator of $W'$ in $A$. As $W'$ is 
a simple submodule of $M$, Lemma 4.7(ii) yields $P\ds_HQ'$ and $Q'\ds_HP$. 
Since $W'\n=0$, we get $Q'\cap Z=\n$.

Suppose now that (b) holds. Then $U\cong(U')^*$ for some $U'\in\M^H\!$,
$\,\dim U'<\infty$. By Lemma 4.6 $M\cong U'\ot V$, so $M$ is finitely 
generated in $\M_A$ according to \cite{Sk-Oy07, Lemma 1.1}. The direct summand 
$M_\n$ of $M$ is also finitely generated in $\M_A$. Then $M_\n$ has a simple 
factor module, call it $W'$. We now complete the proof similarly to case (a), 
but using Lemma 4.7(iii) instead.
\endproof

We say that $A$ has {\it quasilocal localizations} with respect to a central 
subring $Z$ if $A_\m$ is quasilocal and $\m A_\m\sbs\Jac(A_\m)$ for each
$\m\in\Max Z$.

\proclaim
Proposition 4.9.
Suppose that $A$ has quasilocal localizations with respect to $Z$. If either 
all rings $A_\m/\m A_\m$, $\,\m\in\Max Z$, are primary or the antipode of $H$ 
is bijective then the relation $\ds_H$ is symmetric on $\Max A$.
\endproclaim

\Proof.
Since $A_\n$ is quasilocal for any $\n\in\Max Z$, there is a single ideal in 
$\Max_\n A$. Hence $Q'=Q$ in the notation of Lemma 4.8.
\endproof

\proclaim
Corollary 4.10.
If $A$ is commutative then $\ds_H$ is symmetric on $\Max A$.
\endproclaim

\Proof.
The hypotheses of Proposition 4.9 are satisfied if we take $Z=A$.
\endproof

\section
5. Projectivity result for module algebras

Let $A$ be an $H\!$-module algebra and $M\in\AH$. The compatibility of the two 
module structures on $M$ is expressed as
$$
h(va)=\sum\,\,(h\1v)(h\2a)\qquad
{\rm for}\ h\in H,\ v\in M,\ a\in A.
$$

\proclaim
Lemma 5.1.
Suppose that $C$ is a subcoalgebra of $H$ and $I$ an ideal of $A$ such that 
$A/I$ is weakly finite. If $M/MI\cong(A/I)^n$ in $\M_A$ and the $A$-module 
$M/MI_C$ is $n$-generated, then $M/MI_C\cong(A/I_C)^n$ in $\M_A$.
\endproclaim

\Proof.
We will regard $\Hom(C,M/MI)$ as a right $\Hom(C,A/I)$-module by means of the 
convolution action. If $\xi:C\to A/I$ and $\eta:C\to M/MI$ are linear maps, then
$$
(\eta\xi)(c)=\sum\,\,\eta(c\1)\xi(c\2),\qquad c\in C.
$$
For each ideal $J$ of $A$ denote by $\pi_J:M\to M/MJ$ the canonical 
projection. Define $\mh,\mt\in\Hom(C,M/MI)$ for each $m\in M$ by the rules
$$
\mh(c)=\ep(c)\pi_I(m),\qquad\mt(c)=\pi_I(cm)
$$
Pick $e_1,\ldots,e_n\in M$ such that $\pi_I(e_1),\ldots,\pi_I(e_n)$ are a 
basis for the $A/I$-module $M/MI$. Given $\eta\in\Hom(C,M/MI)$, there 
are uniquely determined $\xi_1,\ldots,\xi_n\in\Hom(C,A/I)$ such that 
$\eta(c)=\sum_{i=1}^n\pi_I(e_i)\xi_i(c)$ for all $c\in C$, which is equivalent 
to $\eta=\sum_{i=1}^n\eh_i\xi_i$. Hence $\eh_1,\ldots,\eh_n$ are a basis for 
the $\Hom(C,A/I)$-module $\Hom(C,M/MI)$.

Pick any elements $v_1,\ldots,v_n$ generating $M$ modulo $MI_C$. Given
$m\in M$, there exist $\ze_1,\ldots,\ze_n\in\Hom(H,A)$ such that
$$ 
S(h)m\equiv\sum_{i=1}^n\,v_i\ze_i(h)\mod{MI_C}
$$
for all $h\in H$. Note that $C(MI_C)\sbs MI=\Ker\pi_I$. Taking $c\in C$, we get
$$
\eqalign{
\ep(c)m=\sum_{(c)}\,c_{(1)}S(c_{(2)}\!)m
&{}\equiv\sum_{(c)}\,c_{(1)}\bigl(\sum_{i=1}^nv_i\ze_i(c_{(2)}\!)\bigr)\cr
&{}\equiv\sum_{i=1}^n\sum_{(c)}\,\,(c_{(1)}v_i)\bigl(c_{(2)}\ze_i(c_{(3)}\!)\bigr)
\mod{MI},
}
$$
and applying $\pi_I$, we deduce
$\mh(c)=\sum_{i=1}^n\sum_{(c)}\,\vt_i(c_{(1)})\th_i(c_{(2)})$ where the 
$\th_i$'s are linear maps $C\to A/I$ defined by the formula
$\th_i(c)=\sum_{(c)}c_{(1)}\ze_i(c_{(2)}\!)+I$. This shows that
$\mh=\sum_{i=1}^n\vt_i\,\th_i$. In particular, the submodule of $\Hom(C,M/MI)$ 
generated by $\vt_1,\ldots,\vt_n$ contains $\eh_1,\ldots,\eh_n$. So
$\vt_1,\ldots,\vt_n$ generate the whole $\Hom(C,M/MI)$. Since 
the algebra $\Hom(C,A/I)$ is weakly finite \cite{Sk07, Lemma 7.1}, 
$\vt_1,\ldots,\vt_n$ are in fact a basis for $\Hom(C,M/MI)$ over $\Hom(C,A/I)$.

Suppose that $x_1,\ldots,x_n\in A$ are any elements such that 
$\sum_{i=1}^nv_ix_i\in MI_C$. Then
$$
\sum_{i=1}^n\sum_{(c)}\,(c_{(1)}v_i)(c_{(2)}x_i)
=c\bigl(\,\sum_{i=1}^nv_ix_i\bigr)\in C(MI_C)\sbs MI
$$
for all $c\in C$. Applying $\pi_I$, we rewrite this as 
$\sum_{i=1}^n\vt_i\xt_i=0$ where $\xt_i\in\Hom(C,A/I)$ is defined as in 
section 3, i.e., $\xt_i(c)=cx_i+I$ for $c\in C$. We must have $\xt_i=0$, 
i.e., $x_i\in I_C$ for each $i=1,\ldots,n$. Hence 
$\pi_{I_C}(v_1),\ldots,\pi_{I_C}(v_n)$ are linearly independent over $A/I_C$.  
\endproof

Further on we assume that $A$ has semilocal localizations with respect to a 
central subring $Z$. For a nonnegative $r\in\Q$ and a finitely generated right 
$A$-module $M$ the open subsets $U_r(M)\sbs\Max Z$ and the ideals $I_r(M)$ of 
$A$ were defined in section 1. When $M\in\HM_A$, we use the same notation 
ignoring the $H$-module structure. Denote by $J_r(M)$ the smallest $H$-stable 
ideal of $A$ containing $I_r(M)$.

\setitemsize(ii)
\proclaim
Lemma 5.2.
Let $r=r_P(M)$ where $P\in\Max A$ and $M\in\HM_A$ is an $A$-finite object. 
Suppose that $r_Q(M)\le r$ for each $Q\in\Max A$ such that $P\ds_HQ$. Let 
$r=n/l$ for some integers $l>0$, $\,n\ge0$. Then{\rm:}

\item(i)
$r_Q(M)=r$ for each $Q\in\Max A$ with $P\ds_HQ$.

\item(ii)
$(M/MP_C)^l\cong(A/P_C)^n$ in $\M_A$ for subcoalgebras $C$ of $H$
with $\dim C<\infty$.

\endproclaim

\Proof.
The ring $A/P$ is simple artinian by Lemma 2.1. Then $(M/MP)^l\cong(A/P)^n$ in 
$\M_A$ since the two $A/P$-modules here have equal lengths. Let $C$ be given 
as in (ii). According to Lemma 3.6 $A/P_C$ is a semilocal ring. Since 
$r_Q(M)l\le rl=n$ for any $Q\in\Max A$ with $P_C\sbs Q$, it follows from Lemma 
1.6(i) that the $A/P_C$-module $(M/MP_C)^l$ is $n$-generated. Now Lemma 5.1 
applied to $M^l\in\HM_A$ yields the isomorphism in (ii). Then 
$(M/MQ)^l\cong(A/Q)^n$, and so $r_Q(M)=r$, for any $Q\in\Max A$ with $P_C\sbs 
Q$. As this holds for all finite dimensional subcoalgebras $C$, we deduce (i).  
\endproof

\proclaim
Lemma 5.3.
Let $r=r_P(M)$ where $P\in\Max A$ and $M\in\HM_A$ is an $A$-finite object. 
Suppose that $U_r(M)$ is quasicompact and $r_Q(M)\le r$ for each $Q\in\Max A$ 
such that $Q\cap Z=Q'\cap Z$ for some $Q'\in\Max A$ with $P\ds_HQ'$. Then 
$J_r(M)\sbs P$.
\endproclaim

\Proof.
The isomorphism in Lemma 5.2(ii) enables us to apply Lemma 2.5(ii) with $R=A$
and $K=P_C$. We conclude that $I_r(M)\sbs P_C$ for each finite dimensional 
subcoalgebra $C$ of $H$, whence $I_r(M)\sbs P_H$. Since $P_H$ is an $H$-stable 
ideal of $A$, it follows that $J_r(M)\sbs P_H\sbs P$.
\endproof

Recall from Lemma 2.4 that $r(M)=\sup\{r_P(M)\mid P\in\Max A\}$.

\proclaim
Proposition 5.4.
Given any $A$-finite object $M\in\HM_A$ and $P\in\Max A$ one has $r_P(M)=r(M)$ 
if and only if\/ $P\sps J_{r(M)}(M)$. Moreover, $J_{r(M)}(M)\ne A$.  
\endproclaim

\Proof.
Let $r=r(M)$. We have $r_Q(M)\le r$ for all $Q\in\Max A$. So by Lemma 2.3
$U_r(M)=\Max Z$, which is a quasicompact space. If $r_P(M)=r$ then $P$ 
satisfies the hypothesis of Lemma 5.3, whence $J_r(M)\sbs P$. Conversely, if 
$J_r(M)\sbs P$, then also $I_r(M)\sbs P$, whence $r_P(M)=r$ by Lemma 2.5(iii) 
(where we take $R=A$, $K=P$). Lemma 2.4 says that $r_P(M)=r$ for at least one 
$P\in\Max A$. Hence $J_r(M)\ne A$.
\endproof

\setitemsize(iii)
\proclaim
Corollary 5.5.
Suppose that $M\in\HM_A$ is an $A$-finite object and $A$ has a maximal ideal 
$P$ such that $r_P(M)=r(M)$ and $P$ contains no nonzero $H$-stable ideals of 
$A$. Let $r(M)=n/l$ for some integers $n\ge0$, $l>0$. Then{\rm:}

\item(i)
$r_Q(M)=r(M)$ for all $Q\in\Max A$.

\item(ii)
$M$ is projective in $\M_A;$
$M$ is a generator in $\M_A$ provided $M\ne0$.

\item(iii)
$M^l_\p\cong A_\p^n$ in $\M_{A_\p}$ for each $\p\in\Spec Z$.

\endproclaim

\Proof.
Let $r=r(M)$. Since $J_r(M)$ is an $H$-stable ideal of $A$ contained in $P$, 
we get $J_r(M)=0$; so $I_r(M)=0$ too. Then $J_r(M)\sbs Q$ for any
$Q\in\Max A$, whence (i) holds by Proposition 5.4. Since $U_r(M)=\Max Z$,
we have $M^l_\m\cong A_\m^n$ for any $\m\in\Max Z$ by Proposition 1.11(ii).
If $\p\in\Spec Z$, then $\p\sbs\m$ for some $\m\in\Max Z$. Since $A_\p$ is a 
localization of $A_\m$, we have $M_\p\cong M_\m\ot_{A_\m}A_\p$, whence (iii).
Note that $M$ is projective or a generator in $\M_A$ if and only if so is 
$M^l$. Hence Lemma 2.6 applied to $M^l$ establishes (ii).
\endproof

\proclaim
Theorem 5.6.
Suppose that $A$ is an $H$-simple $H$-module algebra which has semilocal 
localizations with respect to a central subring $Z$. Let $M$ be any locally 
$A$-finite object of $\HM_A$. Put $l=\gcd\{\lng A/Q\mid Q\in\Max A\}$. 
Then{\rm:}

\item(i)
$M$ is projective in $\M_A;$
$M$ is a generator in $\M_A$ provided $M\ne0$.

\item(ii)
$M^l_\p$ is a free $A_\p$-module for each $\p\in\Spec Z$.

\item(iii)
If $M$ is not $A$-finite then $M_\p$ is a free $A_\p$-module for each 
$\p\in\Spec Z$.

\endproclaim

\Proof.
If $M$ is $A$-finite, then there exists $P\in\Max A$ with $r_P(M)=r(M)$ by 
Lemma 2.4. We may now apply Corollary 5.5. For each $Q\in\Max A$ we have
$r(M)=r_Q(M)$, whence $r(M)\cdot\lng(A/Q)\in\Z$. It follows that $r(M)l\in\Z$,
and so Corollary 5.5 establishes both (i) and (ii).  

Suppose further that $M$ is not $A$-finite. The family $\F$ of all 
$\HM_A$-subobjects of $M$ clearly satisfies condition (a) of Lemma 5.7 below. 
If $N\in\F$ and $N\ne M$ then, since $M$ is locally $A$-finite, there exists a 
nonzero $A$-finite subobject $F\sbs M$ such that $F\not\sbs N$. We have
$N'=N+F\in\F$ and $N$ is properly contained in $N'$. Furthermore,
$N'/N\cong F/(F\cap N)$ is an $A$-finite object of $\HM_A$; as we have proved 
already, $N'/N$ is projective in $\M_A$. Thus condition (b) of Lemma 5.7 is 
also fulfilled, and (i) follows.

Let $\m\in\Max Z$. If $N\in\F$ and $N\ne M$ then $M/N$ is a generator in 
$\M_A$ by (i). In this case $(M/N)_\m$ is a generator in $\M_{A_\m}$, whence
$N_\m\ne M_\m$. In particular, this holds for any $A$-finite subobject of $M$ 
since $M$ is not $A$-finite. As a consequence, the $A_\m$-module $M_\m$ cannot 
be finitely generated. The freeness of $M_\m$ now follows from Lemma 5.8 which 
we apply by considering the family of submodules $N_\m$ of $M_\m$ with 
$N\in\F$. If $\p\in\Spec Z$, then $M_\p\cong M_\m\ot_{A_\m}A_\p$ for any 
$\m\in\Max Z$ containing $\p$. This proves (iii).
\endproof

\setitemsize(a)
\proclaim
Lemma 5.7.
Let $R$ be any ring. A right $R$-module $M$ has to be projective provided 
that there exists a family $\F$ of submodules of $M$ satisfying

\item(a)
$\{0\}\in\F$ and the union of every chain in $\F$ is again in $\F${\rm,}

\item(b)
each $N\in\F$, $\,N\ne M$, is properly contained in some $N'\in\F$ such that
$N'/N$ is projective in $\M_R$.

\noindent
If at least one $N\in\F$ is a generator of $\M_R$ then $M$ is a generator too.
\endproclaim

\Proof.
Let $\xi:V\to W$ be any epimorphism and $\ph:M\to W$ any morphism in $\M_R$. 
By Zorn's Lemma there exist a maximal element in the set $X$ of all pairs
$(N,\psi)$ where $N\in\F$ and $\psi:N\to V$ is an $\M_R$-morphism 
such that $\xi\circ\psi=\ph|_N$. If $N$, $N'$ are as in (b), then
$N'=N\oplus G$ for some projective submodule $G$; it is then clear that
any $\psi$ occurring as a component of $(N,\psi)\in X$ can be extended to an 
$\M_R$-morphism $\psi':N'\to V$ with the property that $(N',\psi')\in X$.  
Therefore every maximal element of $X$ has to be $(M,\psi)$ where $\psi:M\to 
V$ is an $\M_R$-morphism satisfying $\xi\circ\psi=\ph$. This proves that $M$ 
is projective. Moreover, the $R$-module $M/N$ is projective for each $N\in\F$ 
since the family of submodules $N'/N$ with $N'\in\F$ and $N'\sps N$ satisfies 
(a) and (b). Hence each $N\in\F$ is a direct summand of $M$, and the final 
assertion of the lemma is clear.
\endproof

\proclaim
Lemma 5.8.
Let $R$ be a semilocal ring. A right $R$-module $M$ is necessarily free as long
as $M$ is not finitely generated and there is a family $\F$ of submodules 
satisfying

\item(a)
$\{0\}\in\F$ and the union of every chain in $\F$ is again in $\F${\rm,}

\item(b)
each $N\in\F\!$, $N\ne M$, is properly contained in some $N'\in\F$ such
that $(N'/N)^l$ is a finitely generated free $R$-module for some $l\in\Z_+$.

\endproclaim

This is a restatement of \cite{Sk07, Lemma 2.5}.

\proclaim
Lemma 5.9.
Let $M\in\HM_A$ be an $A$-finite object. Suppose that $M\ne0$ and $A$ is not 
$H$-simple. Then $A$ has a nonzero $H$-stable ideal $I$ such that $MI\ne M$.  
\endproclaim

\Proof.
Suppose that $MI=M$ for each nonzero $H$-stable ideal $I$ of $A$. Since 
$M\ne0$, we have $r(M)>0$. If $P$ is any maximal ideal of $A$ for which 
$r_P(M)=r(M)$, then $MP\ne M$, and therefore $P$ cannot contain nonzero 
$H$-stable ideals of $A$. Now Corollary 5.5 shows that $M$ is a generator in 
$\M_A$. Then $MI\ne M$ for each proper ideal $I$ of $A$. It follows that 
$A$ cannot have $H$-stable ideals other than $0$ and $A$, i.e. $A$ is 
$H$-simple.
\endproof

\section
6. Local projectivity and flatness

Here we consider an $H$-module algebra $A$ which is not $H$-simple, but there 
is a prime ideal of $A$ containing no nonzero $H$-stable ideals. We want to 
look at the localizations $M_\p$ at a single prime of $Z$. In contrast to 
Theorem 5.6 we are able to prove the projectivity of $M_\p$ only under 
additional restrictions. 

\setitemsize(iii)
\proclaim
Proposition 6.1.
Let $A$ be an $H$-module algebra which has semilocal localizations with 
respect to $Z$. Suppose that $\Max Z$ is noetherian and either all rings 
$A_\m/\m A_\m$, $\,\m\in\Max Z$, are semiprimary or the antipode of $H$ is 
bijective. Let $M\in\HM_A$ be an $A$-finite object whose rank function 
$Q\mapsto r_Q(M)$ is constant on each fibre $\Max_\m A$, $\,\m\in\Max Z$. 
Let $r\in\Q$ and $P\in\Max A$. Then{\rm:}

\item(i)
$r_P(M)=r$ if and only if $P\sps J_r(M)$ and $P\not\sps J_s(M)$ for any $s>r$.

\item(ii)
$r_Q(M)=r_P(M)$ for each $Q\in\Max A$ satisfying $P\ds_HQ$.

Assuming that $P$ contains no nonzero $H$-stable ideals of $A$ and 
$r_P(M)=n/l$ for some integers $n\ge0$, $l>0$, we also have{\rm:}

\item(iii)
$r_Q(M)\ge r_P(M)$ for all $Q\in\Max A$.

\item(iv)
$M_\n^l\cong A_\n^n$ in $\M_{A_\n}$ for any $\n\in\Max Z$ such that 
$r_Q(M)=r_P(M)$ on $\Max_\n A$.

\endproclaim

\Proof.
Since $\Max Z$ is noetherian, for any real $x>0$ the open subset 
$\bigcup_{s<x}U_s(M)$ is quasicompact. Hence there exists $t\in\Q$, $t<x$,
such that $U_s(M)=U_t(M)$ for each $s\in\Q$ satisfying $t<s<x$. Given 
$Q\in\Max A$ and $\n=Q\cap Z$, we have, by Lemma 2.3, $\n\in U_s(M)$ if and 
only if $r_Q(M)\le s$ since the rank function of $M$ is constant on
$\Max_\n A$. It follows that $r_Q(M)\le t$ whenever $r_Q(M)<x$.

The previous argument shows that for any subset $X\sbs\Max A$ there exists
$P'\in X$ such that $r_Q(M)\le r_{P'}(M)$ for all $Q\in X$. For, if we let 
$x=\sup\{r_Q(M)\mid Q\in X\}$ and take $P'$ with $r_{P'}(M)$ sufficiently 
close to $x$, we must have $r_{P'}(M)=x$.

Now choose $P'$ as above in the subset $X=\{Q\in\Max A\mid P\ds_HQ\}$. Denote
$x=r_{P'}(M)$. We have $P\ds_H P'$. If $Q\in\Max A$ satisfies $P'\ds_H Q$, 
then also $P\ds_H Q$, i.e., $Q\in X$. By the assumption on the rank function 
of $M$ we get $r_{Q'}(M)=r_Q(M)\le x$ for any $Q'\in\Max A$ with $Q'\cap 
Z=Q\cap Z$. Thus $P'$ satisfies the hypotheses of Lemmas 5.2, 5.3. We deduce 
that $J_x(M)\sbs P'$ and $r_Q(M)=x$ for any $Q\in\Max A$ with $P'\ds_H Q$.

By Lemma 4.8 there exists $P''\in\Max A$ such that $P''\cap Z=P'\cap Z$, while 
both $P\ds_H P''$ and $P''\ds_H P$ hold. The first condition on $P''$ shows 
that $r_{P''}(M)=x$, while the second condition gives $P''\in X$. But then we 
may replace $P'$ with $P''$ and conclude that $r_Q(M)=x$ for any $Q\in\Max A$
with $P''\ds_H Q$. In particular, $r_P(M)=x$. Now we may replace $P'$ with $P$.
The earlier conclusions about $P'$ yield (ii) and verify the inclusion 
$J_x(M)\sbs P$.

If $s\in\Q$ is such that $r_P(M)<s$ then $r_Q(M)<s$ for all $Q\in\Max_\m A$
where $\m=P\cap Z$. Lemma 2.5(iii) applied with $R=A$, $K=P$ shows that 
$I_s(M)\not\sbs P$; then also $J_s(M)\not\sbs P$ for such $s$. But we have 
checked already that $J_s(M)\sbs P$ for $s=r_P(M)$. The last two statements 
are equivalent to (i).

Suppose that $P$ contains no nonzero $H$-stable ideals of $A$. Then 
$J_x(M)=0$. By (i) applied to an arbitrary $Q\in\Max A$, the inclusion 
$J_x(M)\sbs Q$ yields $r_Q(M)\ge x$, proving (iii). If $\n$ is as in (iv),
then $\n\in U_x(M)$. Since $I_x(M)=0$, Lemma 1.11(ii) verifies (iv).
\endproof

\proclaim
Proposition 6.2.
Let $A$ be an $H$-module algebra, module-finite over a central subring $Z$
such that $Z$ is a Jacobson ring with a noetherian space $\Max Z$. Suppose 
that $P\in\Spec A$ contains no nonzero $H$-stable ideals of $A$. Let 
$M\in\HM_A$ be an $A$-finite object whose rank function $Q\mapsto r_Q(M)$ is 
constant on each fibre $\Max_\m A$, $\,\m\in\Max Z$. Put $m=\inf\{r_Q(M)\mid 
Q\in\Max A\},$ and let $m=n/l$ for some integers $n\ge0$, $l>0$. Then there 
exists $z\in Z,$ $\,z\notin P,$ such that{\rm:}

\item(i)
$r_Q(M)=m$ for each $Q\in\Max A$ with $z\notin Q$.

\item(ii)
$M^l_z\cong A^n_z$ in $\M_{A_z}${\rm;} hence
$M^l_\q\cong A_\q^n$ in $\M_{A_\q}$ for each $\q\in\Spec Z$ with $z\notin\q$.

\endproclaim

\Proof.
Let $\p=P\cap Z$; clearly $\p\in\Spec Z$. As we pointed out in the Remark at 
the end of section 2, $A_\p$ is semilocal. For $r\in\Q$ we have 
$\p\in\Ut_r(M)$ if and only if $r_{Q'}(M_\p)\le r$ for all $Q'\in\Max A_\p$ 
(this follows from Lemma 1.6 and the definition of $\Ut_r(M)$). Hence there 
exists the smallest $r$ with the previous property, namely 
$r=\max\{r_{Q'}(M_\p)\mid Q'\in\Max A_\p\}$. We will assume that $r$ is this 
number.

Since $\Ut_r(M)$ is an open neighborhood of $\p$ in $\Spec Z$, there exists
a basic open subset $D(z)\sbs\Ut_r(M)$ for some $z\in Z\setm\p$. We choose 
such a $z$. If $Q$ is a maximal ideal of $A$ with $z\notin Q$, then 
$Q\in\Max_\n A$ where $\n\in D(z)\cap\Max Z\sbs U_r(M)$, so that $r_Q(M)\le r$.

By Lemma 2.8 there exists a subset $X\sbs\Max A$ such that 
$P=\bigcap_{Q\in X}Q$ and $z\notin Q$ for each $Q\in X$. Suppose that $Q\in X$ 
and $\n=Q\cap Z$. Since $P\sbs Q$, we have $\p\sbs\n$. If $s=r_Q(M)$, then 
$\n\in U_s(M)$, and so $\Ut_s(M)$ is an open neighbourhood of $\n$ in
$\Spec Z$. It follows that $\p\in\Ut_s(M)$, which yields $s\ge r$ by the 
choice of $r$. Since the opposite inequality has been established, we conclude 
that $r_Q(M)=r$. Now $J_r(M)\sbs Q$ by Proposition 6.1(i).

It follows that $J_r(M)\sbs P$. Since $P$ contains no nonzero $H$-stable 
ideals, we get $J_r(M)=0$. Then $I_r(M)=0$ too. Thus $J_r(M)\sbs Q$ for any 
$Q\in\Max A$; Proposition 6.1(i) ensures that $r_Q(M)\ge r$. Note that 
$U_r(M)\ne\varnothing$ since $X\ne\varnothing$. It follows that $m=r$, and the 
previous inequalities prove (i).

Lemma 1.11(ii) shows that $M_\n^l\cong A_\n^n$ when $\n\in U_r(M)$.
Since $Z$ is a Jacobson ring and its localization $Z_z$ at $z$ is a finitely 
generated $Z$-algebra, every maximal ideal of $Z_z$ contracts to a maximal 
ideal of $Z$ \cite{Bou, Ch.~V, \S3, Th.~3}. Thus the maximal ideals of $Z_z$ 
are of the form $\n Z_z$ with $\n\in\Max Z$, $z\notin\n$. We know that
$\n\in U_r(M)$ for any such $\n$. It follows that $A_z$ has semilocal 
localizations with respect to $Z_z$ and the right $A_z$-modules $M^l_z$ and 
$A^n_z$ have isomorphic localizations at all maximal ideals of $Z_z$. Thus we 
may apply Lemma 2.6 to the $A_z$-module $M^l_z$. Replacing $z$ with a 
suitable element $z'$ such that $Z_{z'}$ is a localization of $Z_z$, we prove 
(ii).
\endproof

A restriction on the rank function is a serious deficiency of Propositions 
6.1, 6.2. This restriction is void in the case where all sets $\Max_\m A$ are
singletons.

\proclaim
Theorem 6.3.
Let $A$ be an $H$-module algebra, module-finite over a central subring $Z$
such that $Z$ is a Jacobson ring with a noetherian space $\Max Z$ and each
maximal ideal of $Z$ is contained in a single maximal ideal of $A$. Suppose 
that $P\in\Spec A$ contains no nonzero $H$-stable ideals of $A$. Denote
$\p=P\cap Z$ and
$$
l=\gcd\{\lng A_\p/Q'\mid Q'\in\Max A_\p\}.
$$
Then $M^l_\p$ is a free $A_\p$-module for any locally $A$-finite $M\in\HM_A$. 
\endproclaim

\Proof.
When $M$ is $A$-finite, we may apply Proposition 6.2. The $A_\p$-module 
$M_\p^{l'}$ is free for some $l'\in\Z_+$. Hence $r=r_{Q'}(M_\p)$ does not 
depend on $Q'\in\Max A_\p$. Since $A_\p$ is semilocal, $M_\p^{l'}$ is free in 
$\M_{A_\p}$ for any $l'\in\Z_+$ such that $rl'\in\Z$. Since $r\cdot\lng 
A_\p/Q'\in\Z$ for any $Q'\in\Max A_\p$, we have $rl\in\Z$, whence the 
conclusion.

Suppose that $M$ is not $A$-finite. If $N$, $N'$ are any two $\HM_A$-subobjects
of $M$ such that $N'/N$ is $A$-finite then $(N'_\p/N_\p)^l$ is a free 
$A_\p$-module. If the $A_\p$-module $M_\p$ is finitely generated, then 
$M_\p=N_\p$ for some $A$-finite subobject, and the conclusion is clear. 
Otherwise we apply Lemma 5.8 by considering the family of submodules $N_\p$ of 
$M_\p$ with $N$ running through all $\HM_A$-subobjects of $M$.
\endproof

\proclaim
Theorem 6.4.
Let $B$ be any $H$-module algebra, $A$ an $H$-stable subalgebra contained in 
the center of $B$. Suppose that $A$ is a Jacobson ring with a noetherian space 
$\Max A$ and $IB=B$ for each nonzero $H$-stable ideal $I$ of $A$. Then each 
locally $A$-finite object $M\in\HM_B$ is flat in $\M_A$.
\endproclaim

\Proof.
Given a monomorphism $\ph:V\to W$ in $\M_A$, denote by $K$ the kernel of the
map $\id\ot\ph:M\ot_AV\to M\ot_AW$. Since the latter map may be regarded as an
$\M_B$-morphism, $K$ is a $B$-module. Suppose that $x\in K$ is a nonzero 
element. Denote by $\a$ the annihilator of $x$ in $A$. Then $x$ is annihilated 
by the ideal $\a B$ of $B$, and therefore $\a B\ne B$. There exists $Q\in\Max 
B$ such that $\a B\sbs Q$. Now $\p=Q\cap A$ is a prime ideal of $A$ and 
$\a\sbs\p$. Since $\p B\ne B$, none of the nonzero $H$-stable ideals of $A$ 
can be contained in $\p$. Since $A$ is commutative, we may apply Theorem 6.3
with $Z=A$ and $P=\p$. We deduce that $M_\p$ is projective in $\M_{A_\p}$, 
which implies that the map
$$
\id\ot\ph\ot\id:M\ot_AV\ot_AA_\p\to M\ot_AW\ot_AA_\p
$$
is injective. On the other hand, the kernel of this map coincides with 
$K\ot_AA_\p$ since $A_\p$ is flat in $\M_A$. Thus $K\ot_AA_\p=0$. Then $x$ is 
annihilated by an element in $Z\setm\p$, i.e. $\a\not\sbs\p$. This 
contradiction shows that $K=0$.
\endproof

\section
7. Dualization to comodule algebras

Let $H$ be a bialgebra and $A$ a right $H$-comodule algebra. An object of 
$\MAH$ will be called {\it$A$-finite} if it is finitely generated in $\M_A$. 
An arbitrary object $M$ is a directed union of its $A$-finite subobjects. 
Indeed, any finite subset of $M$ is contained in a finite dimensional 
$H$-subcomodule; the $A$-submodule generated by the latter is an $A$-finite 
subobject.

\proclaim
Lemma 7.1.
Each object of $\MAH$ is flat {\rm(}resp. projective{\rm)} in $\M_A$ provided 
that this is true for all $A$-finite objects. Each nonzero object of $\MAH$ is 
a projective generator in $\M_A$ provided that this is true for all nonzero
$A$-finite objects.
\endproclaim

\Proof.
Since tensor products commute with filtered direct limits, the flat part of 
the lemma follows from the fact that each $M\in\MAH$ is a directed union of 
$A$-finite subobjects. The projective part follows from Lemma 5.7 in which we 
take $\F$ to be the family of all subobjects of $M$.
\endproof

Let $H'$ be a second bialgebra, $A'$ an $H'$-comodule algebra. Given 
a homomorphism of bialgebras $\ph:H'\to H$, we may view $A'$ as an 
$H$-comodule algebra. Suppose that we are given also a map $A'\to A$ which is 
a homomorphism of $H$-comodule algebras. In this case there is a functor 
$\M_{A'}^{H'}\rightsquigarrow\MAH$ which takes an object $N\in\M_{A'}^{H'}$ to 
$N\ot_{A'}A\in\MAH$ on which the comodule structure 
$N\ot_{A'}A\to(N\ot_{A'}A)\ot H$ is given by the rule
$$
v\ot a\mapsto\sum\,(v\0\ot a\0)\ot\ph(v\1)a\1
$$
where $v\in N$, $a\in A$. It is easy to check that this map is well-defined.
In the special case where $H'=H$ and $A'=k$ with the trivial comodule 
structure, we obtain an object $V\ot A\in\MAH$ for each right $H$-comodule 
$V$.

We next make several observations concerning direct limits of comodule 
algebras. Suppose that $H=\limdir H_i$, the direct limit of an inductive 
family $\H=(H_i)$ of bialgebras indexed by a directed set $\I$. An
{\it$\H$-compatible inductive family of comodule algebras} $\F=(A_i)$ is a 
collection containing for each $i\in\I$ an $H_i$-comodule algebra $A_i$ and 
for each pair $i,j\in\I$ with $i\le j$ a homomorphism of $H_j$-comodule 
algebras $A_i\to A_j$; these maps are requested to obey the usual rules of
inductive systems. If such an $\F$ is given, $A=\limdir A_i$ becomes an 
$H$-comodule algebra in a natural way. We mention below several properties
of the category $\MAH$ under previous assumptions.

We say that $M\in\MAH$ is {\it$\F$-induced} if there exists $i\in\I$ and 
$N\in\M_{A_i}^{H_i}$ such that $M\cong N\ot_{A_i}A$. Denote by $\FAH$ the 
class of all $A$-finite objects of $\MAH$ isomorphic to $N\ot_{A_i}A$ for some 
$i\in\I$ and an $A_i$-finite $N\in\M_{A_i}^{H_i}$.

\proclaim
Lemma 7.2.
If $M\in\MAH$ is an $\F$-induced object and $M'$ any $A$-finite subobject,
then $M/M'$ is $\F$-induced. In this case $M/M'\in\FAH$ whenever $M\in\FAH$.
\endproclaim

\Proof.
Let $M\cong N\ot_{A_i}A$ for some $i\in\I$ and $N\in\M_{A_i}^H$. Put
$J=\{j\in\I\mid i\le j\}$. For each $j\in J$ denote $N_j=N\ot_{A_i}A_j$, and 
let $M_j$ be the image of the canonical map $\ph_j:N_j\to M$. We thus obtain a 
directed family of vector subspaces of $M$ indexed by $J$. Since $A$ is 
covered by the images of $A_j$, $j\in J$, we have $M=\bigcup_{j\in J}M_j$. By 
the hypothesis $M'$ is generated in $\M_A$ by a finite subset. The latter is 
contained in a finite dimensional $H$-subcomodule $V\sbs M'$, and we then have 
$M'=VA$. There exists $j\in J$ such that $V\sbs M_j$. We may view $\ph_j$ as 
an $\M^H$-morphism. Hence there exists a finite dimensional $H$-subcomodule $W$
of $N_j$ such that $\ph_j(W)=V$. Let $\rho_j:N_j\to N_j\ot H_j$ be the
$H_j$-comodule structure map. We must have
$$
\rho_j(W)\sbs W\ot H_j+N_j\ot\Ker(H_j\to H).
$$
Then $\rho_j(W)\sbs W\ot H_j+N_j\ot U$ for some finite dimensional subspace in 
the kernel of $H_j\to H$. Now $U$ vanishes in $H_t$ for some $t\in J$,
$\,t\ge j$. If we denote by $W'$ the image of $W$ in $N_t$, then $W'$ is
an $H_t$-subcomodule of $N_t$ satisfying $\ph_t(W')=V$. The map 
$W'\ot A_t\to N_t$ afforded by the $A_t$-module structure is a morphism in 
$\M_{A_t}^{H_t}$; hence its cokernel $K$ is an object of that category. 
Tensoring with $A$, we obtain an exact sequence
$$
(W'\ot A_t)\ot_{A_t}A\to N_t\ot_{A_t}A\cong M\to K\ot_{A_t}A\to0
$$
in $\MAH$. By construction the image of the first map coincides with $M'$.
It follows that $M/M'\cong K\ot_{A_t}A$ is an $\F$-induced object. Note that 
$K$ is $A_t$-finite whenever $N$ is $A_i$-finite.
\endproof

\proclaim
Lemma 7.3.
Every $A$-finite object $M\in\MAH$ is isomorphic to a factor object of an 
object from $\FAH$.
\endproclaim

\Proof.
There exists a finite dimensional $H$-subcomodule $V\sbs M$ such that $M=VA$.
The map $V\ot A\to M$ afforded by the $A$-module structure is then an 
epimorphism in $\MAH$. So it remains to prove the conclusion of the lemma for
the object $V\ot A$.

Since $H$ is an injective cogenerator in $\M^H$, we can embed $V$ as a
subcomodule in $H^n$ for some integer $n>0$. For each $i\in\I$ let 
$\ph_i:H_i^n\to H^n$ denote the canonical map. Then $V$ is contained in the 
image of $\ph_j$ for some $j\in\I$. Since $\ph_j$ may be regarded as an 
$\M^H$-morphism, there exists a finite dimensional $H$-subcomodule $W\sbs 
H_j^n$ such that $\ph_j(W)=V$. As in the proof of previous lemma we can find 
$t\in\I$, $\,t\ge j$, such that the image $W'$ of $W$ in $H_t^n$ is an 
$H_t$-subcomodule. The map $W'\to V$ obtained by restriction of $\ph_t$ is an 
epimorphism in $\M^H$; it gives rise to an epimorphism $W'\ot A\to V\ot A$ in 
$\MAH$. Since $W'\ot A\cong(W'\ot A_t)\ot_{A_t}A$ and $W'\ot A_t$ is an 
$A_t$-finite object of $\M_{A_t}^{H_t}$, we have $W'\ot A\in\FAH$.
\endproof

\proclaim
Lemma 7.4.
Suppose $A$ has a homomorphism into a nonzero artinian ring $R$ and each 
nonzero object from $\FAH$ is a projective generator in $\M_A$. Then $\FAH$ 
contains all $A$-finite objects of $\MAH$.
\endproclaim

\Proof.
For each $N\in\FAH$ the $R$-module $N\ot_AR$ is finitely generated, and we 
denote by $l(N)$ its length. If $N\ne0$, then $N$ is a projective generator in 
$\M_A$; in this case $N\ot_AR$ is a projective generator in $\M_R$, whence
$N\ot_AR\ne0$, i.e. $l(N)>0$.

If $M\in\MAH$ is $A$-finite, Lemma 7.3 ensures the existence of an 
$\MAH$-epimorphism $\ph:F\to M$ with $F\in\FAH$. If $M$ is projective in 
$\M_A$, then $\ph$ splits in $\M_A$, and therefore $\Ker\ph$ is an $A$-finite 
object of $\MAH$. In this case $M\in\FAH$ by Lemma 7.2.

Suppose now that $M\in\FAH$ and $M'$ is any $A$-finite $\MAH$-subobject of $M$.
By Lemma 7.2 $M/M'\in\FAH$. By the hypothesis $M$ and $M/M'$ are both 
projective in $\M_A$. Then so is $M'$ too. This implies $M'\in\FAH$ as we have 
observed above. Since the exact sequence $0\to M'\to M\to M/M'\to0$ splits 
in $\M_A$, it remains exact after tensoring with $R$, whence 
$l(M)=l(M')+l(M/M')$. If $M'\ne M$ then $l(M/M')>0$, in which case 
$l(M')<l(M)$.

Given another $A$-finite subobject $M''$ of $M$ properly containing $M'$, we 
have then $l(M')<l(M'')$ since $M''$ is also in $\FAH$. It is now clear that 
$M$ satisfies ACC on $A$-finite subobjects. But every subobject of 
$M$ is a directed union of $A$-finite ones; hence it is itself $A$-finite.
So, according to Lemma 7.2, the class $\FAH$ is closed under factor 
objects, and we are done.
\endproof

For the ring $R$ appearing in the next lemma we say that a finitely generated 
projective $R$-module $G$ has {\it constant rank} if $r_P(G)$, as defined in 
section 1, does not depend on $P\in\Max R$; we denote by $r(G)$ this common 
value. If $G\ne0$, then $G$ has a simple factor module annihilated by some 
$P$, and therefore $r(G)=r_P(G)>0$. If $G\cong G'\oplus G''$ in $\M_R$, then
$r_P(G)=r_P(G')+r_P(G'')$ for all $P$; hence $G''$ has constant rank whenever 
so do both $G$ and $G'$. In this case $r(G')<r(G)$ unless $G''=0$. There are 
only finitely many possible values of $r(G')$ when $G'$ runs through the 
direct summands of $G$ having constant rank; indeed, $r(G')<r(G)$ and 
$r(G')l\in\Z$ for any $G'$ where $l$ is the greatest common divisor of the 
lengths of the simple artinian factor rings of $R$. It follows that $G$ 
satisfies ACC on direct summands of constant rank.

\proclaim
Lemma 7.5.
Let $A\to R$ be a homomorphism into a ring $R$ all whose right primitive 
factor rings are artinian. Suppose that for each $F\in\FAH$ the $R$-module 
$F\ot_AR$ is projective of constant rank. Then for each $A$-finite $N\in\MAH$ 
there exists an epimorphism $\xi:F\to N$ in $\MAH$ such that $F\in\FAH$ and 
the map $F\ot_AR\to N\ot_AR$ induced by $\xi$ is an isomorphism in $\M_R$.
\endproclaim

\Proof.
Let $M\in\FAH$. By the hypothesis $M\ot_AR$ is a projective $R$-module of 
constant rank. For each $\MAH$-subobject $K\sbs M$ denote by $T_K$ the image
of the canonical map $K\ot_AR\to M\ot_AR$. If $M'\sbs M$ is an $A$-finite 
subobject, then $M/M'\in\FAH$ by Lemma 7.2. In this case $M/M'\ot_AR$ is a
projective $R$-module of constant rank, and it follows from the exact sequence
$$
M'\ot_AR\to M\ot_AR\to M/M'\ot_AR\to0
$$
that so too is $T_{M'}$. As a consequence, $M\ot_AR$ satisfies ACC on 
submodules of the form $T_{M'}$ with $M'$ as above. An arbitrary subobject
$K\sbs M$ is a directed union of $A$-finite ones. Then $T_K=\bigcup\,T_{K'}$ 
where $K'$ runs through the $A$-finite subobjects of $K$, and it follows that
$T_K=T_{K'}$ for some $K'$ of this type. Then the canonical projection
$\xi:M/K'\to M/K$ induces an isomorphism after tensoring with $R$. Thus for 
$N=M/K$ we have the desired conclusion with $F=M/K'$. Lemma 7.3 completes the 
proof.
\endproof

\setitemsize(iii)
\proclaim
Proposition 7.6.
Let $H=\limdir H_i$ and $A=\limdir A_i$ as before.

\item(i)
All objects of $\MAH$ are flat in $\M_A$ provided that for each $i$
all objects of $\M_{A_i}^{H_i}$ are flat in $\M_{A_i}$.

\item(ii)
Suppose that $A$ has a homomorphism into a nonzero artinian ring $R$.
If for each $i$ all nonzero objects of $\M_{A_i}^{H_i}$ are projective 
generators in $\M_{A_i}$ then all nonzero objects of $\MAH$ are projective 
generators in $\M_A$.

\item(iii)
Let $A\to R$ be a ring homomorphism where $R$ is a ring all whose right
primitive factor rings are artinian. Suppose that for each $i$ and each 
$A_i$-finite $N\in\M_{A_i}^{H_i}$ the $R$-module $N\ot_{A_i}R$ is projective 
of constant rank. Then $M\ot_AR$ is projective in $\M_R$ for any 
$M\in\MAH${\rm;} if $M$ is $A$-finite then $M\ot_AR$ has constant rank.

\endproclaim

\Proof.
(i) Since $N\ot_{A_i}A$ is flat in $\M_A$ whenever $N$ is flat in $\M_{A_i}$, 
the hypothesis implies that every $\F$-induced object of $\MAH$ is flat in 
$\M_A$. An arbitrary $A$-finite object $L\in\MAH$ is isomorphic to $M/K$ where 
$M\in\FAH$ and $K$ is an $\MAH$-subobject of $M$. We have $L\cong\limdir M/K'$ 
where $K'$ ranges over all $A$-finite subobjects of $K$. By Lemma 7.2
$M/K'\in\FAH$ for each $K'$. Thus $L$ is a direct limit of flat $A$-modules; 
then $L$ is flat in $\M_A$. Lemma 7.1 completes the proof.

(ii) The hypothesis implies that all nonzero $\F$-induced objects of $\MAH$ 
are projective generators in $\M_A$. Lemma 7.4 shows that the same conclusion 
holds for all $A$-finite objects, and Lemma 7.1 establishes this for 
arbitrary objects.

(iii) Here $M\ot_AR$ is a projective $R$-module of constant rank for each
$M\in\FAH$. By Lemma 7.5 the same holds for each $A$-finite object of 
$\MAH$. In order to extend this to arbitrary objects we have to repeat the 
proof of \cite{Sk07, Th. 1.2} given in case of commutative algebras. The proof 
is easier under the assumption that $R$ is flat in $\AM$. In this case we can 
apply Lemma 5.7 to the family of submodules $N\ot_AR$ of $M\ot_AR$ where $N$ 
runs through all $\MAH$-subobjects of $M$.
\endproof

Recall that the finite dual $\Hd$ of $H$ is a subalgebra of $H^*$ consisting 
of all linear functions vanishing on an ideal of finite codimension in $H$.
There is a comultiplication on $\Hd$ dual to the multiplication on $H$. 
Moreover, $\Hd$ is a Hopf algebra whenever so is $H$. As explained in 
\cite{Sw}, $\M^H$ is equivalent to the category of rational left $H^*$-modules.
This gives a functor $\M^H\rightsquigarrow\HdM$. If $A$ is a right $H$-comodule
algebra, then $A$ is a left $\Hd$-module algebra with respect to the 
corresponding module structure; then we obtain a functor 
$\MAH\rightsquigarrow\HdMA$. Moreover, all objects in the image of that 
functor are locally $A$-finite.

If $H$ is residually finite dimensional, then $\Hd$ is dense in $H^*$; it 
follows that the subcomodules of any $U\in\M^H$ coincide with the submodules 
of the corresponding $\Hd$-module. In this case the $H$-costable ideals of an 
$H$-comodule algebra $A$ coincide with the $\Hd$-stable ideals, and $A$ is an 
$H$-simple $H$-comodule algebra if and only if $A$ is an $\Hd$-simple 
$\Hd$-module algebra.

It is easy now to translate the results from the preceding sections to the 
context of comodule algebras. The next result is the comodule version of 
Theorem 5.6.

\proclaim
Theorem 7.7.
Let $M\in\MAH$ where $H$ is a residually finite dimensional Hopf algebra and 
$A$ is an $H$-simple $H$-comodule algebra which has semilocal localizations 
with respect to a central subring $Z$. Put $l=\gcd\{\lng A/Q\mid Q\in\Max A\}$.
Then{\rm:}

\item(i)
$M$ is projective in $\M_A;$
$M$ is a generator in $\M_A$ provided $M\ne0$.

\item(ii)
$M^l_\p$ is a free $A_\p$-module for each $\p\in\Spec Z$.

\item(iii)
If $M$ is not $A$-finite then $M_\p$ is a free $A_\p$-module for each 
$\p\in\Spec Z$.

\endproclaim

In case of commutative comodule algebras we can weaken the assumption about 
the Hopf algebra $H$.

\proclaim
Theorem 7.8.
Suppose $H$ is a directed union of residually finite dimensional Hopf 
subalgebras and $A$ is a commutative $H$-comodule algebra. If $\p\in\Spec A$ 
contains no nonzero $H$-costable ideals of $A$ then $M_\p$ is a free 
$A_\p$-module for any $M\in\MAH$.
\endproclaim

\Proof.
If $H$ is residually finite dimensional and $A$ is finitely generated, the 
conclusion follows from Theorem 6.3 since in this case $Z=A$ is a noetherian 
Jacobson ring. In general there is a directed family $\G$ of residually finite 
dimensional Hopf subalgebras of $H$ whose union coincides with $H$. Let $\I$ 
be the set of all pairs $(A',H')$ where $H'\in\G$ and $A'$ is a finitely 
generated subalgebra of $A$ such that $\rho_A(A')\sbs A'\ot H'$. Given two 
pairs $(A_1,H_1)$ and $(A_2,H_2)$ from $\I$, there exists $H_3\in\G$ 
containing both $H_1$ and $H_2$; clearly $(A_1A_2,H_3)\in\I$. This shows that 
$\I$ is directed by inclusion.

For $(A',H')\in\I$ we may regard $A'$ as an $H'$-comodule algebra. Any
$H'$-costable ideal of $A'$ extends to an $H$-costable ideal of $A$. It 
follows then that the prime ideal $\p'=\p\cap A'$ of $A'$ contains no nonzero 
$H'$-costable ideals of $A'$. Hence $N\ot_{A'}A'_{\p'}$ is a free 
$A'_{\p'}$-module for any $N\in\M_{A'}^{H'}$; since the homomorphism
$A'\to A_\p$ factors through $A'_{\p'}$, the $A_\p$-module $N\ot_{A'}A_\p$ is 
also free. The assignment $(A',H')\mapsto H'$ defines an inductive family 
$\H$ of Hopf algebras indexed by $\I$. The direct limit of $\H$ is equal to 
$H$. The assignment $(A',H')\mapsto A'$ defines an $\H$-compatible inductive 
family $\F$ of comodule algebras. If $V$ is any $H$-subcomodule of $A$ with
$\dim V<\infty$, then $\rho_A(V)\sbs V\ot H'$ for some $H'\in\G$; denoting by 
$A'$ the subalgebra of $A$ generated by $V$, we have $(A',H')\in\I$. Since $A$
is covered by its finite dimensional subcomodules, the direct limit of $\F$ 
equals $A$.

It remains to apply Proposition 7.6(iii) with $R=A_\p$. Since the ring $A_\p$ 
is local, all projective $A_\p$-modules are free by Kaplansky's Theorem.
\endproof

Theorem 7.8 implies the next result whose proof follows that of Theorem 6.4.

\proclaim
Theorem 7.9.
Let $H$ be a directed union of residually finite dimensional Hopf subalgebras, 
$B$ an $H$-comodule algebra, and $A$ an $H$-costable subalgebra contained in 
the center of $B$. Suppose that $IB=B$ for each nonzero $H$-costable ideal $I$ 
of $A$. Then each object $M\in\MBH$ is flat in $\M_A$.
\endproclaim

In conclusion we will deduce all results stated in the introduction. The 
structure theorem for objects of $\M^H_H$ \cite{Sw, Th. 4.1.1} shows that any 
Hopf algebra $H$ is a simple object of $\M^H_H$. If $A$ is a Hopf subalgebra 
of $H$, then $A$ is simple in $\M_A^A$, and therefore simple in $\MAH$. In 
this case $A$ is an $H$-simple $H$-comodule algebra.

If $A$ is a right coideal subalgebra of $H$, then $IH=H$ for each $H$-costable 
ideal $I\ne\nobreak0$ of $A$ since $IH$ is an $\M^H_H$-subobject of $H$.  
Furthermore, the opposite multiplication makes $A\op$ into a right coideal 
subalgebra of the bialgebra $H\op$. If the antipode of $H$ is bijective, 
$H\op$ is actually a Hopf algebra. The previous argument applied to $A\op$, 
$H\op$ yields $HI=H$ for each nonzero $H$-costable ideal $I$ of $A$.

\proclaim
Theorem 7.10.
Let $H$ be a residually finite dimensional Hopf algebra, $A$ be a Hopf 
subalgebra having semilocal localizations with respect to a central 
subring $Z$. Then each nonzero object $M\in\MAH$ is a projective generator in 
$\M_A$ and $M_\p$ is a free $A_\p$-module for any $\p\in\Spec Z$.
\endproclaim

\Proof.
We apply Theorem 7.7 in which $l=1$ since $k$ is a factor algebra of $A$.
\endproof

In particular, $H$ is a projective generator in $\M_A$. We may change both the 
multiplication and comultiplication in $A$ and $H$ to the opposite ones, 
obtaining another pair of Hopf algebras $A^{\rm op,cop}\sbs H^{\rm op,cop}$. 
Theorem 7.10 applied to the latter shows that $H$ is a projective generator in 
$\AM$. Thus Theorem 0.1 is proved.

\proclaim
Theorem 7.11.
Let $A\sbs B\sbs H$ where $H$ is a residually finite dimensional Hopf algebra,
$B$ is a Hopf subalgebra, and $A$ is a right coideal subalgebra which has 
semilocal localizations with respect to a central subring $Z$. Suppose that 
$B$ is right module-finite over $A$ and the antipode of $B$ is bijective. Then 
each nonzero object $M\in\MAH$ is a projective generator in $\M_A$ and $M_\p$ 
is a free $A_\p$-module for any $\p\in\Spec Z$.
\endproclaim

\Proof.
We have $BI=B$ for each nonzero $H$-costable ideal $I$ of $A$. Hence we may 
apply Lemma 5.9 regarding $M=B$ as an $A$-finite object of $\HdMA$. It follows 
that $A$ is an $H$-simple $H$-comodule algebra. Again Theorem 7.7 applies.  
\endproof

\proclaim
Theorem 7.12.
Let $A\sbs B\sbs H$ where $H$ is a directed union of residually 
finite dimensional Hopf subalgebras, $B$ is any Hopf subalgebra, and $A$ is a 
right coideal subalgebra contained in the center of $B$. Then{\rm:}

\item(i)
$M_\p$ is a free $A_\p$-module for each $M\in\MAH$ and $\p\in\Spec A$ 
with $\p B\ne B$.

\item(ii)
Each object of $\MBH$ is flat in $\M_A$.

\endproclaim

\Proof.
Since $\De(A)\sbs(A\ot H)\cap(B\ot B)=A\ot B$, we may regard $A$ as a right 
coideal subalgebra of $B$. Furthermore, an ideal $I$ of $A$ is $H$-costable if 
and only if it is $B$-costable. It follows that $IB=B$ for each nonzero 
$H$-costable ideal $I$ of $A$. In particular, $\p\in\Spec A$ cannot contain
nonzero $H$-costable ideals of $A$ whenever $\p B\ne B$. The two conclusions 
are therefore consequences of Theorems 7.8, 7.9.
\endproof

\proclaim
Theorem 7.13.
Let $A$ be a commutative Hopf subalgebra of a Hopf algebra $H$ which is a 
directed union of residually finite dimensional Hopf subalgebras.
Then each nonzero object of $\MAH$ is a projective generator in $\M_A$.
\endproclaim

\Proof.
If $H$ is residually finite dimensional then the conclusion is a special case 
of Theorem 7.10. In general we apply this to each residually finite 
dimensional Hopf subalgebra $H'$ of $H$ and the right coideal subalgebra 
$A\cap H'$ of $H'$. Proposition 7.6(ii) completes the proof.
\endproof

\Remark.
If $A$ is contained in the center of $H$ then Theorem 7.13 can be proved by 
first observing that $H$ is faithfully flat in $\M_A$ and then using 
\cite{Tak79, Th. 5}.
\endremark

\references
\nextref
Az51
\auth
G.,Azumaya;
\endauth
\paper{On maximally central algebras}
\journal{Nagoya Math.~J.}
\Vol{2}
\Year{1951}
\Pages{119-150}

\nextref
Bass
\auth
H.,Bass;
\endauth
\book{Algebraic $K$-theory}
\publisher{Benjamin}
\Year{1968}

\nextref
Bj73
\auth
J.-E.,Bj\"ork;
\endauth
\paper{Noetherian and Artinian chain conditions of associative rings}
\journal{Arch. Math.}
\Vol{24}
\Year{1973}
\Pages{366-378}

\nextref
Bou
\auth
N.,Bourbaki;
\endauth
\book{Commutative Algebra}
\publisher{Springer}
\Year{1989}

\nextref
Br
\auth
K.A.,Brown;K.R.,Goodearl;
\endauth
\book{Lectures on Algebraic Quantum Groups}
\publisher{Birkh\"auser}
\Year{2002}

\nextref
Camps93
\auth
R.,Camps;W.,Dicks;
\endauth
\paper{On semilocal rings}
\journal{Isr. J. Math.}
\Vol{81}
\Year{1993}
\Pages{203-211}

\nextref
Ch90
\auth
W.,Chin;
\endauth
\paper{Spectra of smash products}
\journal{Isr. J. Math.}
\Vol{72}
\Year{1990}
\Pages{84-98}

\nextref
Cor88
\auth
B.,Cortzen;L.W.,Small;
\endauth
\paper{Finite extensions of rings}
\journal{Proc. Amer. Math. Soc.}
\Vol{103}
\Year{1988}
\Pages{1058-1062}

\nextref
Cur53
\auth
C.W.,Curtis;
\endauth
\paper{Noncommutative extensions of Hilbert rings}
\journal{Proc. Amer. Math. Soc.}
\Vol{4}
\Year{1953}
\Pages{945-955}

\nextref
Con94
\auth
C.,De Concini;V.,Lyubashenko;
\endauth
\paper{Quantum function algebra at roots of $1$}
\journal{Adv. Math.}
\Vol{108}
\Year{1994}
\Pages{205-262}

\nextref
Doi83
\auth
Y.,Doi;
\endauth
\paper{On the structure of relative Hopf modules}
\journal{Comm. Algebra}
\Vol{11}
\Year{1983}
\Pages{243-255}

\nextref
Et04
\auth
P.,Etingof;V.,Ostrik;
\endauth
\paper{Finite tensor categories}
\journal{Moscow Math.~J.}
\Vol{4}
\Year{2004}
\Pages{627-654}

\nextref
Hoff
\auth
K.,Hoffmann;
\endauth
\paper{Coidealunteralgebren in endlich dimensionalen Hopfalgebren}
Dissertation, Univ. M\"unchen, 1991.

\nextref
Kh
\auth
V.K.,Kharchenko;
\endauth
\paper{PBW-bases of coideal subalgebras and a freeness theorem}
Preprint.

\nextref
Lam
\auth
T.-Y.,Lam;
\endauth
\book{A First Course in Noncommutative Rings}
\bookseries{Graduate Texts in Mathematics}
\Vol{131}
\publisher{Springer}
\Year{1991}

\nextref
Ma91
\auth
A.,Masuoka;
\endauth
\paper{On Hopf algebras with cocommutative coradicals}
\journal{J. Algebra}
\Vol{144}
\Year{1991}
\Pages{451-466}

\nextref
Ma92
\auth
A.,Masuoka;
\endauth
\paper{Freeness of Hopf algebras over coideal subalgebras}
\journal{Comm. Algebra}
\Vol{20}
\Year{1992}
\Pages{1353-1373}

\nextref
Ma94
\auth
A.,Masuoka;D.,Wigner;
\endauth
\paper{Faithful flatness of Hopf algebras}
\journal{J. Algebra}
\Vol{170}
\Year{1994}
\Pages{156-164}

\nextref
Mo
\auth
S.,Montgomery;
\endauth
\book{Hopf algebras and Their Actions on Rings}
\bookseries{CBMS Regional Conference Series in Mathematics}
\Vol{82}
\publisher{American Mathematical Society}
\Year{1993}

\nextref
Mo95
\auth
S.,Montgomery;H.-J.,Schneider;
\endauth
\paper{Hopf crossed products, rings of quotients, and prime ideals}
\journal{Adv. Math.}
\Vol{112}
\Year{1995}
\Pages{1-55}

\nextref
Mo99
\auth
S.,Montgomery;H.-J.,Schneider;
\endauth
\paper{Prime ideals in Hopf Galois extensions}
\journal{Isr. J.~Math.}
\Vol{112}
\Year{1999}
\Pages{187-235}

\nextref
Ni89
\auth
W.D.,Nichols;M.B.,Zoeller;
\endauth
\paper{A Hopf algebra freeness theorem}
\journal{Amer. J. Math.}
\Vol{111}
\Year{1989}
\Pages{381-385}

\nextref
Pro67
\auth
C.,Procesi;
\endauth
\paper{Non commutative Jacobson-rings}
\journal{Ann. Scuola Norm. Sup. Pisa}
\Vol{21}
\Year{1967}
\Pages{281-290}

\nextref
Rad77
\auth
D.E.,Radford;
\endauth
\paper{Pointed Hopf algebras are free over Hopf subalgebras}
\journal{J. Algebra}
\Vol{45}
\Year{1977}
\Pages{266-273}

\nextref
Rob81
\auth
J.C.,Robson;L.W.,Small;
\endauth
\paper{Liberal extensions}
\journal{Proc. London Math. Soc.}
\Vol{42}
\Year{1981}
\Pages{87-103}

\nextref
Scha00
\auth
P.,Schauenburg;
\endauth
\paper{Faithful flatness over Hopf subalgebras: counterexamples}
\In{Interactions between ring theory and representations of algebras}
\bookseries{Lect. Notes Pure Appl. Math.}
\Vol{210}
\publisher{Marcel Dekker}
\Year{2000}
\Pages{331-344}

\nextref
Schn93
\auth
H.-J.,Schneider;
\endauth
\paper{Some remarks on exact sequences of quantum groups}
\journal{Comm. Algebra}
\Vol{21}
\Year{1993}
\Pages{3337-3357}

\nextref
Sk07
\auth
S.,Skryabin;
\endauth
\paper{Projectivity and freeness over comodule algebras}
to appear in
\journal{Trans. Amer. Math. Soc.}

\nextref
Sk-Oy07
\auth
S.,Skryabin;F.,Van Oystaeyen;
\endauth
\paper{The Goldie Theorem for H-semiprime algebras}
\journal{J.~Algebra}
\Vol{305}
\Year{2006}
\Pages{292-320}

\nextref
Sw
\auth
M.E.,Sweedler;
\endauth
\book{Hopf Algebras}
\publisher{Benjamin}
\Year{1969}

\nextref
Tak71
\auth
M.,Takeuchi;
\endauth
\paper{Free Hopf algebras generated by coalgebras}
\journal{J.~Math. Soc. Japan}
\Vol{23}
\Year{1971}
\Pages{561-582}

\nextref
Tak79
\auth
M.,Takeuchi;
\endauth
\paper{Relative Hopf modules---equivalences and freeness criteria}
\journal{J. Algebra}
\Vol{60}
\Year{1979}
\Pages{452-471}

\nextref
Wu03
\auth
Q.-S.,Wu;J.J.,Zhang;
\endauth
\paper{Noetherian PI Hopf algebras are Gorenstein}
\journal{Trans. Amer. Math. Soc.}
\Vol{355}
\Year{2003}
\Pages{1043-1066}

\endreferences
\bye